\documentclass[12pt]{article}
\usepackage{amsmath,amssymb,amsfonts,amsthm,url,xspace,graphicx,psfrag}
\usepackage[pdftitle={Combinatorics of Tripartite Boundary Connections for Trees and Dimers},
            pdfauthor={Richard W. Kenyon and David B. Wilson},
            bookmarks=true,bookmarksopen=true,bookmarksopenlevel=2]{hyperref}
\usepackage[all]{hypcap}

\newcommand{\sref}[1]{\hyperref[#1]{Section~\ref*{#1}}}
\newcommand{\aref}[1]{\hyperref[#1]{Appendix~\ref*{#1}}}
\newcommand{\lref}[1]{\hyperref[#1]{Lemma~\ref*{#1}}}
\newcommand{\tref}[1]{\hyperref[#1]{Theorem~\ref*{#1}}}
\newcommand{\cref}[1]{\hyperref[#1]{Corollary~\ref*{#1}}}
\newcommand{\fref}[1]{\hyperref[#1]{Figure~\ref*{#1}}}
\newcommand{\pref}[1]{\hyperref[#1]{Proposition~\ref*{#1}}}
\newcommand{\eref}[1]{\hyperref[#1]{Equation~\ref*{#1}}}

\def\clap#1{\hbox to 0pt{\hss#1\hss}}

 \newcommand{\MRhref}[2]{\href{http://www.ams.org/mathscinet-getitem?mr=#1}{MR#2}}

\makeatletter
\def\@strippedMR{}
\def\@scanforMR#1#2#3\endscan{%
  \ifx#1M\ifx#2R\def\@strippedMR{#3}%
  \else\def\@strippedMR{#1#2#3}%
  \fi\fi}
\def\@rst#1 #2other{#1}
\newcommand\MR[1]{\relax\ifhmode\unskip\spacefactor3000 \space\fi
  \@scanforMR#1\endscan
  \MRhref{\expandafter\@rst\@strippedMR other}{\@strippedMR}}
\newcommand\MRs[1]{\relax\ifhmode\unskip\spacefactor3000 \space\fi
  \@scanforMR#1\endscan
  \MRhref{\@strippedMR}{\@strippedMR}}
\makeatother

\DeclareMathOperator*{\doublesum}{\sum\ \sum}
\newtheorem{theorem}{Theorem}
\numberwithin{theorem}{section}

\newtheorem{lemma}[theorem]{Lemma}
\newtheorem{proposition}[theorem]{Proposition}

\theoremstyle{definition}
\theoremstyle{definition}
\newcommand{\g}{\mathsf{t}}
\newcommand{\DD}{\mathsf{DD}}
\newcommand{\WB}{\mathsf{WB}}
\newcommand{\BW}{\mathsf{BW}}
\newcommand{\eqz}{\stackrel{\g}{\equiv}}
\newcommand{\eqt}{\stackrel{2}{\equiv}}
\newcommand{\eqq}{\stackrel{q}{\equiv}}
\newcommand{\SLE}{\operatorname{SLE}}
\newcommand{\LHS}{\text{LHS}}
\newcommand{\RHS}{\text{RHS}}
\newcommand{\Z}{\mathbb{Z}}
\newcommand{\R}{\mathbb{R}}

\newcommand{\M}{\mathcal M}

\newcommand{\G}{\mathcal G}
\newcommand{\K}{\mathcal K}
\renewcommand{\H}{\mathbb{H}}
\newcommand{\I}{\mathbf{I}}
\newcommand{\No}{\mathbf{N}}
\newcommand{\V}{\mathbf{V}}
\newcommand{\N}{\mathbb{N}}
\newcommand{\E}{\mathcal E}
\newcommand{\ZDD}{Z^{\DD}}
\newcommand{\ZBW}{Z^{\BW}}
\newcommand{\ZWB}{Z^{\WB}}

\newcommand{\X}{X}
\newcommand{\ZD}{{Z^{\mathsf{D}}}}
\newcommand{\KL}{L}

\newcommand{\eps}{\varepsilon}
\renewcommand{\Im}{\operatorname{Im}}

\newcommand{\old}[1]{}

\renewcommand{\th}{\ensuremath{^{\text{th}}}\xspace}

\newcommand{\unc}{\operatorname{uncrossing}}
\newcommand{\tree}{\operatorname{tree}}

\newcommand{\CdV}{MR1652692}
\newcommand{\DGG}{MR1462755}
\newcommand{\CIM}{MR1657214}
\newcommand{\Fomin}{MR1837248}

\newcommand{\Kasteleyn}{MR0253689}
\newcommand{\LSW}{MR2044671}
\newcommand{\CS}{MR2097339}
\newcommand{\Spitzer}{MR0388547}

\newcommand{\Pt}{\overline\Pr}
\newcommand{\Pu}{\ddddot\Pr}
\newcommand{\hPr}{\widehat\Pr}
\newcommand{\pt}[1]{\Pt(#1)}
\newcommand{\pu}[1]{\Pu(#1)}
\renewcommand{\SLE}{\text{SLE}}

\newcommand{\fm}{\phantom{-}}

\newcommand{\p}{{\mathcal P}^{(\g)}}
\newcommand{\pd}{{\mathcal P}^{(\DD)}}

\def\rcs $#1: #2 ${\expandafter\def\csname rcs#1\endcsname {#2}}
\rcs $Date: 2008/12/17 07:47:29 $

\begin{document}
\title{\vspace*{-80pt}Boundary Partitions in Trees and Dimers
\footnotetext{2000 \textit{Mathematics Subject Classification.}  60C05, 82B20, 05C05, 05C50.}
\footnotetext{\textit{Key words and phrases.}  Tree, grove, double-dimer model, Gaussian free-field, Dirichlet-to-Neumann matrix, meander, SLE.}
}
\author{Richard W. Kenyon\thanks{University of British Columbia, Vancouver, BC V6T 1Z2, Canada} \and David B. Wilson\thanks{Microsoft Research, Redmond, WA 98052, USA}}
\old{
\author{Richard W. Kenyon}
\address{Richard W. Kenyon\\ Department of Mathematics\\ University of British Columbia\\ Vancouver, B.C. V6T 1Z2, Canada}
\author{David B. Wilson}
\address{David B. Wilson\\ Microsoft Research\\ One Microsoft Way\\ Redmond, WA\ 98052, USA}
}
\date{}
\maketitle

\vspace*{-18pt}
\old{
\begin{abstract}
  We study \textbf{groves} on planar graphs, which are forests in
  which every tree contains one or more of a special set of vertices
  on the outer face, referred to as \textbf{nodes}.  Each grove
  partitions the set of nodes.  When a random grove is selected, we
  show how to compute the various partition probabilities as functions
  of the electrical properties of the graph when viewed as a resistor
  network.  We prove that for any partition $\sigma$, Pr[grove has
  type $\sigma$] / Pr[grove is a tree] is a dyadic-coefficient
  polynomial in the pairwise resistances between the nodes, and
  Pr[grove has type $\sigma$] / Pr[grove has maximal number of trees]
  is an integer-coefficient polynomial in the entries of the
  Dirichlet-to-Neumann matrix.  We give analogous integer-coefficient
  polynomial formulas for the pairings of chains in the double-dimer
  model.  We show that the distribution of pairings of contour
  lines in the Gaussian free field with certain natural boundary conditions is
  identical to the distribution of pairings in the scaling limit of the
  double-dimer model.  These partition probabilities are relevant to
  multichordal $\SLE_2$, $\SLE_4$, and $\SLE_8$.
\end{abstract}
}
\begin{abstract}
  Given a finite planar graph, a grove is a spanning forest in which
  every component tree contains one or more of a specified set of
  vertices (called nodes) on the outer face.  For the uniform measure
  on groves, we compute the probabilities of the different possible
  node connections in a grove.  These probabilities only depend on
  boundary measurements of the graph and not on the actual graph
  structure, i.e., the probabilities can be expressed as functions of
  the pairwise electrical resistances between the nodes, or
  equivalently, as functions of the Dirichlet-to-Neumann operator (or
  response matrix) on the nodes.  These formulae can be likened to
  generalizations (for spanning forests) of Cardy's percolation
  crossing probabilities, and generalize Kirchhoff's formula for the
  electrical resistance.  Remarkably, when appropriately normalized,
  the connection probabilities are in fact integer-coefficient
  polynomials in the matrix entries, where the coefficients have a
  natural algebraic interpretation and can be computed combinatorially.
  A similar phenomenon holds in
  the so-called double-dimer model: connection probabilities of
  boundary nodes are polynomial functions of certain boundary
  measurements, and as formal polynomials, they are specializations of
  the grove polynomials.  Upon taking scaling limits, we show that the
  double-dimer connection probabilities coincide with those of the
  contour lines in the Gaussian free field with certain natural
  boundary conditions.  These results have direct application to
  connection probabilities for multiple-strand $\SLE_2$, $\SLE_8$, and
  $\SLE_4$.
\end{abstract}

\section{Introduction}

\subsection{Grove partitions}

A \textbf{circular planar graph} $\G$ is a finite weighted planar graph with
a set of vertices $\No$ on its outer face numbered $1,\dots,n$ in
counterclockwise order.  The vertices in $\No$ are called \textbf{nodes}, and
the remaining vertices are called \textbf{inner vertices}.
Define a \textbf{grove} to be a spanning acyclic subgraph (a forest) of
$\G$ such that each component tree contains at least one node.
The weight of a grove is the product of the weights of the edges it
contains.  We study random groves where the probability of a grove is
proportional to its weight.

The term \textbf{grove} comes from Carroll and Speyer \cite{\CS} and
Petersen and Speyer \cite{MR2144860} who studied a special case of
(our) groves, with a particular family of underlying graphs.  Since we
are dealing with a natural generalization we chose to keep their
terminology, and refer to the special case they discuss as
Carroll-Speyer groves, which we will discuss further in a subsequent
paper \cite{KW-Pfaffian}.

The connected components of a grove partition the nodes into a planar
(i.e.\ noncrossing)
partition.  For example, when $n=4$, there are 14 planar partitions:
$1234$, $1|234$, $2|134$, $3|124$, $4|123$, $12|34$, $23|14$,
$1|2|34$, $1|3|24$, $1|4|23$, $2|3|14$, $2|4|13$, $3|4|12$, $1|2|3|4$.
There are no groves with the partition $13|24$ because it is not
planar (there is no way to connect nodes $1$ and $3$ and nodes $2$ and $4$
by disjoint paths within a circular planar graph).  For general~$n$, the
number of noncrossing partitions is the $n$\th Catalan number
$C_n = (2n)!/(n!(n+1)!)$ (see \cite[ex.~6.19(pp), pg.~226]{MR1676282}).
\enlargethispage{9pt}

\begin{figure}[htbp]
\centerline{
\includegraphics{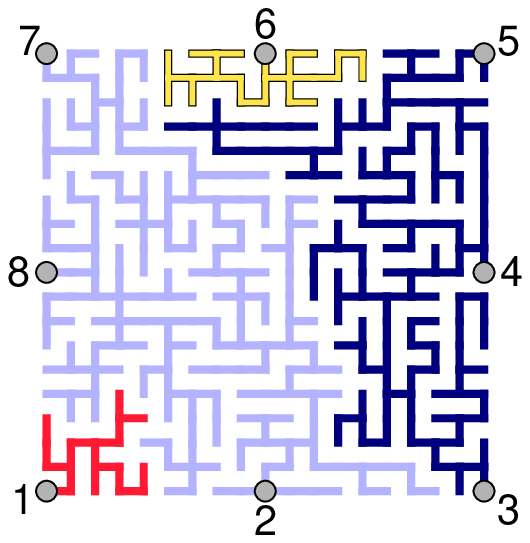}\hfil\includegraphics[viewport=20 -45 11 -29]{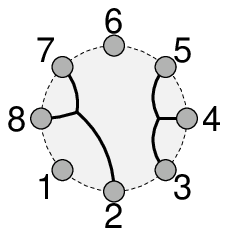}}
\caption{
  A random grove (left) of a rectangular grid with $8$ nodes on the
  outer face.  In this grove there are $4$ trees (each colored
  differently), and the partition of the nodes is $\{\{1\}, \{2,7,8\},
  \{3,4,5\}, \{6\}\}$, which we write as $1|278|345|6$, and illustrate
  schematically as shown on the right.
\label{grove}}
\end{figure}

If $\sigma$ is a planar partition of $1,\dots,n$, we let
$\Pr(\sigma)$ denote the probability that a random grove of $\G$
partitions the nodes according to $\sigma$.  Since groves with one
component are trees, we refer to the partition $\sigma=123\dots n$
as the \textbf{tree} partition.
When there are $n$ nodes we call the partition $\sigma=1|2|3|\cdots |n$
the \textbf{uncrossing}, since groves with this partition type contain
no crossings (i.e.\ paths) connecting the nodes.
We show how to compute for each
planar partition~$\sigma$ the probability $\Pr(\sigma)$ of it occurring
in a random grove, as a function of the electrical properties of the
graph~$\G$, when $\G$ is viewed as a resistor network with
conductances equal to the edge weights.

We let $R_{i,j}$ denote the effective
electrical resistance between nodes $i$ and~$j$, i.e., the voltage
at node~$i$ which, when node~$j$ is held at $0$~volts, causes a unit current
to flow through the circuit from node~$i$ to node~$j$.  (A good reference
for basic electric circuit theory is \cite{MR920811}.)  Let $L_{i,j}$
denote the current that would flow into node~$j$ if node~$i$ were set to
one volt and the remaining nodes set to zero volts.
Though it is not obvious from this definition, $L_{i,j}=L_{j,i}$
(see \aref{response}).
The $L_{i,j}$ are the negatives of the entries of the ``response matrix''
(Dirichlet-to-Neumann matrix) of $(\G,\No)$, on which we provide
further background in \aref{response}, see also \cite{\CdV}.
We prove
\begin{theorem}\label{L-R-thm}
  For any planar partition $\sigma$,
 $$ \frac{\Pr(\sigma)}{\Pr(\tree)}
   = \genfrac{}{}{0pt}{}{\text{integer-coefficient homogeneous polynomial in the $R_{i,j}/2$'s}}{\text{where the degree is $-1+\#\text{parts of $\sigma$}$,}\hfil}$$
and
 $$ \frac{\Pr(\sigma)}{\Pr(\unc)}
   = \genfrac{}{}{0pt}{}{\text{integer-coefficient homogeneous polynomial in the $L_{i,j}$'s}}{\text{where the degree is $n-\#\text{parts of $\sigma$}$.}\hfil}$$
\end{theorem}
This theorem is a corollary to our main theorem on groves,
\tref{P-thm} in \sref{groves-explicit}, which gives
explicit formulas.  These formulas are in effect a multinode
generalization (for planar graphs) of Kirchhoff's formula
$\Pr(1|2)/\Pr(12) = R_{1,2} = 1/L_{1,2}$ (which holds whether or not
$\G$ is planar) \cite{Kirchhoff}
(see \cite[Thm~3.7]{MR1410532} for a good exposition).
  (When $n=2$ we have
$L_{1,2}=1/R_{1,2}$, but this does not hold for $n>2$.)  We illustrate
\tref{L-R-thm} by giving the polynomials for $n=2$ and
$n=3$, and in \aref{4nodes} we give the polynomials for $n=4$
nodes.  For notational convenience we write
\enlargethispage*{1.7\baselineskip}
\begin{align*}
\pt\sigma&:=\Pr(\sigma)/\Pr(\tree)&
&\text{and}& \pu\sigma&:=\Pr(\sigma)/\Pr(\unc).
\end{align*}
(The dots in
$\pu$ remind us that in the normalization the nodes are disconnected.)
The polynomials for $n=2$ come from Kirchhoff's formula:
\begin{align*}
\pt{12}&=1, & \pt{1|2} &= R_{1,2}, & \quad&\text{and}\quad&
\pu{12}&=L_{1,2}, & \pu{1|2} &= 1.
\end{align*}
When $n=3$ the polynomials start to become more interesting:
$$
\begin{matrix}
\pt{123}=1,\ \ \ \ \ \ \ \ 
\pt{1|23} = \frac12 R_{1,3} + \frac12 R_{1,2} - \frac12 R_{2,3},\\
\pt{2|13} = \frac12 R_{2,3} + \frac12 R_{1,2} - \frac12 R_{1,3},\ \ \ \ \ \ \ \ 
\pt{3|12} = \frac12 R_{1,3} + \frac12 R_{2,3} - \frac12 R_{1,2},
\\
\pt{1|2|3}= \frac12 R_{1,2}R_{1,3} +\frac12 R_{1,2}R_{2,3} +\frac12 R_{1,3}R_{2,3}
            -\frac14 R_{1,2}^2 -\frac14 R_{2,3}^2 -\frac14 R_{1,3}^2,
\end{matrix}
$$
and
$$
\begin{matrix}
\pu{123}=L_{1,2}L_{1,3}+L_{1,2}L_{2,3}+L_{1,3}L_{2,3},\\
\pu{1|23} = L_{2,3},\ \ \ \ \ \ 
\pu{2|13} = L_{1,3},\ \ \ \ \ \ 
\pu{3|12} = L_{1,2},\ \ \ \ \ \ 
\pu{1|2|3}= 1.
\end{matrix}
$$
The $n=3$ formulas also hold whether or not the graph is planar,
but when $n\geq 4$ planarity becomes important.

When $n$ gets large these polynomials can have many terms, but for
certain special classes of partitions, most notably the ``parallel
crossing'' $1,n|2,n-1|3,n-2|\cdots$, the polynomials can be
expressed in terms of determinants \cite{\CIM} \cite{\Fomin}.  We
plan to discuss these and other determinant / Pfaffian formulas in
a subsequent article \cite{KW-Pfaffian}.

\subsection{The partition projection matrix and explicit formulas}
\label{groves-explicit}

The formulas in \tref{L-R-thm} in terms of the $L_{i,j}$'s are the simplest ones to explain.
For a partition $\tau$ on $1,\dots,n$ we define
\begin{equation}\label{Ltau}
\KL_\tau = \sum_F\prod_{\text{$\{i,j\} \in F$}} L_{i,j},
\end{equation}
where the sum is over those spanning forests $F$ of the complete graph on
$n$ vertices $1,\dots,n$ for which trees of $F$ span the parts of $\tau$,
and the product is over edges $\{i,j\}$ of forest $F$.

This definition makes sense whether or not the partition $\tau$ is planar.
For example, $\KL_{1|234} = L_{2,3} L_{3,4} + L_{2,3}L_{2,4} + L_{2,4} L_{3,4}$
and $\KL_{13|24} = L_{1,3} L_{2,4}$.
As we shall see, the ``$L$ polynomials'' of \tref{L-R-thm} are in fact
integer linear combinations of the $\KL_\tau$'s:
 $$\frac{\Pr(\sigma)}{\Pr(1|2|\cdots|n)} = \sum_\tau \p_{\sigma,\tau} \KL_\tau.$$
 (We write the superscript $(\g)$ to distinguish these coefficients
 from ones that arise in the double-dimer model in
 \sref{double-dimer-explicit}.)
The rows of the matrix $\p$ are indexed by planar partitions,
and the columns are indexed by all partitions.
In the case of $n=4$ nodes, the matrix $\p$ is

\vspace{3pt}
\newcommand{\rf}[1]{\begin{rotatebox}{90}{$#1$}\end{rotatebox}}
\centerline{
\begin{tabular}{r@{$\,\,\,$}c@{$\,$}c@{$\,$}c@{$\,$}c@{$\,$}c@{$\,$}c@{$\,$}c@{$\,$}c@{$\,$}c@{$\,$}c@{$\,$}c@{$\,$}c@{$\,$}c@{$\,$}c@{$\,$}c}
        & \rf{1|2|3|4} & \rf{12|3|4} & \rf{13|2|4} & \rf{14|2|3} & \rf{23|1|4} & \rf{24|1|3} & \rf{34|1|2} & \rf{12|34} & \rf{14|23} & \rf{1|234} & \rf{2|134} & \rf{3|124} & \rf{4|123} & \rf{1234\phantom{|}} & \rf{13|24} \\
$1|2|3|4$  &1&0&0&0&0&0&0&0&0&0&0&0&0&0& 0 \\
$12|3|4$   &0&1&0&0&0&0&0&0&0&0&0&0&0&0& 0 \\
$13|2|4$   &0&0&1&0&0&0&0&0&0&0&0&0&0&0& 0 \\
$14|2|3$   &0&0&0&1&0&0&0&0&0&0&0&0&0&0& 0 \\
$23|1|4$   &0&0&0&0&1&0&0&0&0&0&0&0&0&0& 0 \\
$24|1|3$   &0&0&0&0&0&1&0&0&0&0&0&0&0&0& 0 \\
$34|1|2$   &0&0&0&0&0&0&1&0&0&0&0&0&0&0& 0 \\
$12|34$    &0&0&0&0&0&0&0&1&0&0&0&0&0&0&$-1\fm$ \\
$14|23$    &0&0&0&0&0&0&0&0&1&0&0&0&0&0&$-1\fm$ \\
$1|234$    &0&0&0&0&0&0&0&0&0&1&0&0&0&0& 1 \\
$2|134$    &0&0&0&0&0&0&0&0&0&0&1&0&0&0& 1 \\
$3|124$    &0&0&0&0&0&0&0&0&0&0&0&1&0&0& 1 \\
$4|123$    &0&0&0&0&0&0&0&0&0&0&0&0&1&0& 1 \\
$1234$     &0&0&0&0&0&0&0&0&0&0&0&0&0&1& 0.
\end{tabular}}
\vspace{3pt}
For example, the row for $1|234$ tells us
$$
  \frac{\Pr(1|234)}{\Pr(1|2|3|4)} = \KL_{1|234} + \KL_{13|24} = L_{2,3} L_{3,4} + L_{2,3}L_{2,4} + L_{2,4} L_{3,4} + L_{1,3}L_{2,4}
$$
and the row for $12|34$ tells us
$$
  \frac{\Pr(12|34)}{\Pr(1|2|3|4)} = \KL_{12|34} - \KL_{13|24} = L_{1,2} L_{3,4} - L_{1,3}L_{2,4}.
$$

We call this matrix $\p$ the \textbf{projection matrix from partitions to
planar partitions}, since it can be interpreted as a map from the
vector space whose basis vectors are indexed by all partitions to the
vector space whose basis vectors are indexed by planar partitions, and
the map is the identity on planar partitions.
For example, the column for $13|24$ tells us
$$
13|24\ \ \text{projects to}\ \ -12|34-14|23+1|234+2|134+3|124+4|123.
$$
(We could have written the right-hand-side as $-e_{12|34}-e_{14|23}+e_{1|234}+e_{2|134}+e_{3|124}+e_{4|123}$, where the $e$'s are basis vectors, but it is convenient to suppress the vector notation and instead write formal linear combinations of partitions.)

The projection matrix may be computed using some simple combinatorial
transformations of partitions.  Given a partition~$\tau$, the
$\tau$\th column of $\p$ may be computed by repeated application of
the following transformation rule, until the resulting formal linear
combination of partitions only involves planar partitions.
The rule generalizes the transformation
\enlargethispage{10pt}
\begin{center}
\begin{psfrags}
\psfrag{+}[cc][Bc][1][0]{$+$}
\psfrag{-}[cc][Bc][1][0]{$-$}
\psfrag{->}[cc][Bc][1][0]{$\to$}
\includegraphics{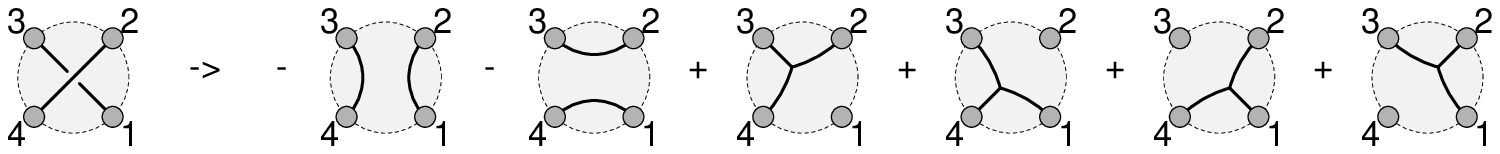}
\end{psfrags}
\end{center}
$$
13|24 \to -12|34-14|23+1|234+2|134+3|124+4|123
$$
(which is derived in \lref{cross})
to partitions $\tau$ containing additional items and parts.
If partition $\tau$ is nonplanar, then there will exist items $a<b<c<d$
such that $a$ and $c$ belong to one part, and $b$ and $d$ belong to
another part.  Arbitrarily subdivide the part containing $a$ and $c$
into two sets $A$ and $C$ such that $a\in A$ and $c\in C$, and
similarly subdivide the part containing $b$ and $d$ into $B\ni b$ and
$D\ni d$.  Let the remaining parts of partition $\tau$ (if any) be
\newcommand{\rest}{\operatorname{rest}}
denoted by ``$\rest$.''  Then the transformation rule is
\begin{multline}\label{part-trans}\tag{Rule~1}
 AC|BD|\rest \to
 A|BCD|\rest
\,+\,B|ACD|\rest
\,+\,C|ABD|\rest
\,+\,D|ABC|\rest\\
-\,AB|CD|\rest
\,-\,AD|BC|\rest.
\end{multline}
In \sref{groves} we prove
\begin{theorem}\label{P-thm}
  Any partition $\tau$ may be transformed into a formal linear
  combination of planar partitions by repeated application of
  \ref{part-trans}, and the resulting linear combination does not
  depend on the choices made when applying \ref{part-trans}, so
  that we may write
  $$\tau \to \sum_{\text{\rm planar partitions $\sigma$}} \p_{\sigma,\tau} \sigma.$$
  For any planar partition $\sigma$,
  these same coefficients $\p_{\sigma,\tau}$ satisfy the equation
 $$\frac{\Pr(\sigma)}{\Pr(1|2|\cdots|n)}
 = \sum_{\text{\rm partitions $\tau$}} \p_{\sigma,\tau} \KL_\tau$$
  for circular planar graphs.
  More generally, for any graph these coefficients satisfy
 $$\sum_{\text{\rm partitions $\tau$}} \p_{\sigma,\tau} \frac{\Pr(\tau)}{\Pr(1|2|\cdots|n)}
 = \sum_{\text{\rm partitions $\tau$}} \p_{\sigma,\tau} \KL_\tau.$$
\end{theorem}
(The last equation specializes to the previous one because $\p_{\sigma,\sigma}=1$, for planar~$\tau\neq\sigma$ we have $\p_{\sigma,\tau}=0$, and for nonplanar~$\tau$ we have $\Pr(\tau)=0$ for circular planar graphs.)

\subsection{Double-dimer pairings}
\label{double-dimer}

\begin{figure}[htbp]
\includegraphics[width=\textwidth]{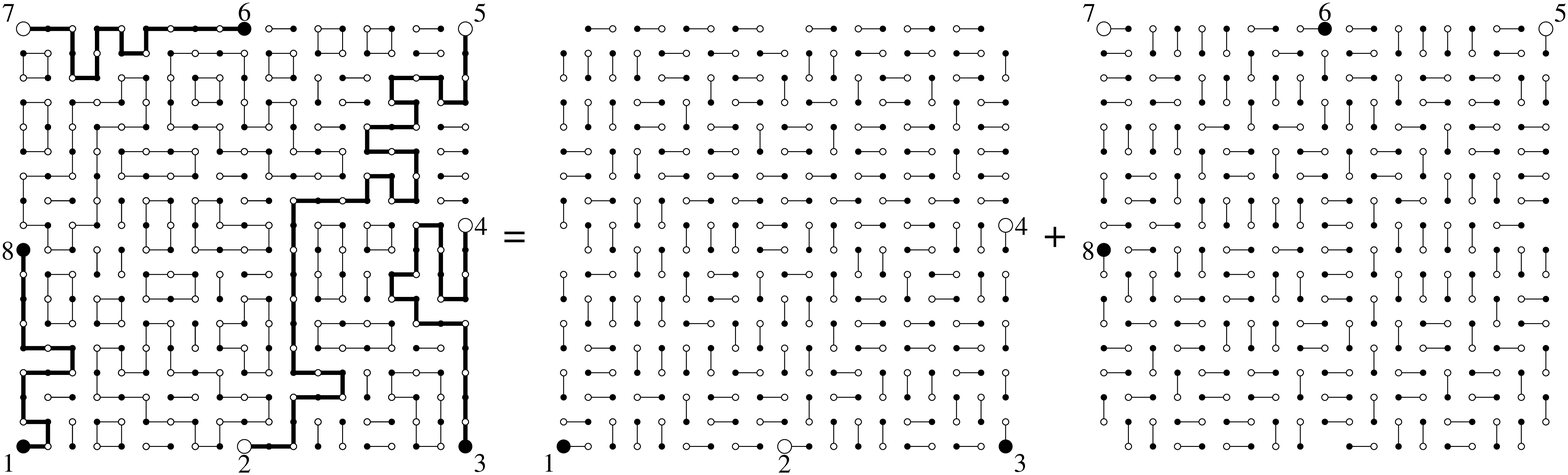}
  \caption{\label{ddimer-dimer}
    At the left is a double-dimer configuration on a rectangular grid
    with $8$ nodes.  In this configuration, the pairing of the
    nodes is $\{\{1,8\},\{3,4\},\{5,2\},\{7,6\}\}$, which we write as
    $\,{}^1_8\!\mid\!{}^3_4\!\mid\!{}^5_2\!\mid\!{}^7_6\,$
    (odd-numbered nodes are always paired with even-numbered nodes).  A
    double-dimer configuration is formed from the union of two dimer
    configurations, one on the graph $\G^\BW\subseteq\G$ (defined below
    and shown in the middle) and the other on the graph $\G^\WB\subseteq \G$
    (shown on the right).
}
\end{figure}

An analogous computation can be done for the double-dimer model.  Let
$\G$ be a finite edge-weighted bipartite planar graph, with a set
$\No$ of $2n$ special vertices called nodes on the outer face, which
we label $1,\dots,2n$ in counterclockwise order along the outer face.
A double-dimer configuration on $(\G,\No)$ is by definition a
configuration of disjoint loops (i.e.\ simple cycles of length more than two),
doubled edges, and simple paths on $\G$ that connect all the nodes in
pairs, which collectively cover the vertices of $\G$.  See
\fref{ddimer-dimer} for an example.  We weight each configuration by
the product of its edge weights times $2^\ell$, where $\ell$ is the
number of loops (a doubled edge does not count as a loop).  The number
of possible ways that the $2n$ nodes may be paired up with one another
is the $n$\th Catalan number (see \cite[ex.~6.19(n), pg.~222]{MR1676282}),
and we show how to compute the
probabilities of each of these pairings when a random double-dimer
configuration is chosen according to these weights.  We will show how
to express these pairing probabilities in terms of polynomials in a
set of variables that are analogous to the $L_{i,j}$ variables for
random groves.

Let $\ZDD(\G,\No)$ be the weighted sum of all double-dimer
configurations.  Since the graph~$\G$ is bipartite, we view the
vertices as being colored black and white so that each edge of $\G$
connects a black and white vertex, and we define the parity of a node
to be the parity of its numerical label.
Let $\G^\BW$ be the subgraph of $\G$ formed by
deleting the nodes except the ones that are black and odd or
white and even, and let $\G^\BW_{i,j}$ be defined as $\G^\BW$ was, but
with node~$i$ included in $\G^\BW_{i,j}$ if and only
if it was not included in $\G^\BW$, and similarly for node~$j$.
Let $\ZBW$ and $\ZBW_{i,j}$ be
the weighted sum of dimer configurations of $\G^\BW$ and
$\G^\BW_{i,j}$ respectively, and define $\ZWB$ and $\ZWB_{i,j}$
similarly but with the roles of black and white reversed.  Each of
these quantities can be computed via determinants, see
\cite{\Kasteleyn} and \sref{ddimers}.
It turns out that $\ZDD = \ZBW \ZWB$; this is essentially Ciucu's
graph factorization theorem \cite[Thm.~1.2]{ciucu} (except that in the factorization theorem,
any edges that connect two nodes are reweighted by $1/2$, and here they are not), see \sref{ddimers}.
(The two dimer configurations in \fref{ddimer-dimer} are on the
graphs $\G^\BW$ and $\G^\WB$.)
The variables that play the role of $L_{i,j}$ in groves are defined by
$$X_{i,j} = \ZBW_{i,j} / \ZBW.$$
\old{
Assuming that the node colors alternate between black
and white, let $\X=\X(\G)$ be the weighted sum of dimer covers of $\G$,
and let $\X_{i,j}=\X(\G\setminus\{i,j\})$ be the weighted sum of dimer
covers of $\G\setminus\{i,j\}$.  Each of $\X$, $\ZDD$, and $\X_{i,j}$ can
be computed via determinants, see \cite{\Kasteleyn} and
\sref{ddimers}.
}
For each planar matching of the nodes
$\sigma$, let $\Pr(\sigma)$ be the probability that a random double-dimer
configuration has this set of connections.  In \sref{ddimers}
we prove
\begin{theorem} \label{Z-thm}
  For any planar pairing $\sigma$ on $2 n$ nodes,
  $$\Pr(\sigma)\frac{\ZWB}{\ZBW} = \text{integer-coefficient homogeneous polynomial of degree $n$ in the quantities $\X_{i,j}$}.$$
\end{theorem}
This theorem is a corollary to our main theorem on double-dimer pairings,
\tref{P2-thm} in \sref{double-dimer-explicit},
which gives the polynomials explicitly.
To illustrate this theorem, we give a few examples.
For notational simplicity let us define $$\hPr(\sigma)=\Pr(\sigma)\ZWB/\ZBW.$$
In the simplest nontrivial case there are $2n=4$ nodes, and two possible planar
pairings, $\{\{1,2\},\{3,4\}\}$ and $\{\{1,4\},\{2,3\}\}$, which
we write as $\,{}^1_2\!\mid\!{}^3_4\,$ and $\,{}^1_4\!\mid\!{}^3_2\,$.
Their normalized probabilities are
\begin{align*}
\hPr(\,{}^1_2\!\mid\!{}^3_4\,)&=\X_{1,2}\X_{3,4} &
\hPr(\,{}^1_4\!\mid\!{}^3_2\,)&=\X_{1,4}\X_{2,3}.
\end{align*}
These above two formulas are essentially equivalent to a formula of Kuo \cite[Thms~2.1 and~2.3]{kuo-condense}.
In the case of $2n=6$ nodes, we have
\begin{align*}
\hPr(\,{}^1_2\!\mid\!{}^3_6\!\mid\!{}^5_4\,)=&\X_{3,6}(\X_{1,2}\X_{4,5}-\X_{1,4}\X_{2,5})\\
\hPr(\,{}^1_2\!\mid\!{}^3_4\!\mid\!{}^5_6)=&\X_{1,4}\X_{2,5}\X_{3,6}+\X_{1,2}\X_{3,4}\X_{5,6}
\end{align*}
and similarly for cyclic permutations of the indices.  These cover the five
possibilities.

In the case of $2n=8$ nodes, we have
\begin{align*}
\hPr(\,{}^1_2\!\mid\!{}^3_8\!\mid\!{}^5_6\!\mid\!{}^7_4\,)=&(\X_{1,2}\X_{7,4}-\X_{1,4}\X_{7,2})(\X_{3,8}\X_{5,6}-\X_{5,8}\X_{3,6})\\
=& -\det\begin{bmatrix}\X_{1,2}&\X_{1,4}&0&0\\
0&0&\X_{3,6}&\X_{3,8}\\0&0&\X_{5,6}&\X_{5,8}\\
\X_{7,2}&\X_{7,4}&0&0\end{bmatrix}\\
\hPr(\,{}^1_2\!\mid\!{}^3_4\!\mid\!{}^5_6\!\mid\!{}^7_8\,)=&
\X_{1,2}\X_{3,4}\X_{5,6}\X_{7,8}+\X_{1,4}\X_{3,8}\X_{5,6}\X_{7,2}+\X_{1,6}\X_{3,4}\X_{5,8}\X_{7,2}+\\
&\X_{1,6}\X_{3,8}\X_{5,2}\X_{7,4}+\X_{1,2}\X_{3,6}\X_{5,8}\X_{7,4}+
\X_{1,4}\X_{3,6}\X_{5,2}\X_{7,8}-2\X_{1,4}\X_{3,6}\X_{5,8}\X_{7,2}
\end{align*}
\begin{align*}
\hPr(\,{}^1_2\!\mid\!{}^3_8\!\mid\!{}^5_4\!\mid\!{}^7_6\,)
=&\det\begin{bmatrix}\X_{1,2}&\X_{1,4}&\X_{1,6}&0\\
0&0&\X_{3,6}&\X_{3,8}\\ \X_{5,2}&\X_{5,4}&0&\X_{5,8}\\
\X_{7,2}&\X_{7,4}&\X_{7,6}&0\end{bmatrix}
\end{align*}
and similarly for cyclic permutations of the indices.  These cover the
fourteen possibilities.

As can be seen from the way we have written these polynomials, for
some pairings they are expressible as determinants.  We plan to discuss
such formulas further in \cite{KW-Pfaffian}.

Of course we could also express the pairing probabilities in terms of
the variables $X^*_{i,j} = \ZWB_{i,j}/\ZWB$.  The polynomials are
exactly the same, although the underlying variables represent different
quantities.

\subsection{The odd-even pairing projection matrix and explicit formulas}
\label{double-dimer-explicit}

The explicit computation of the double-dimer pairing probability
formulas is quite analogous to the computation of the grove partition
probability formulas.
The role of the $\KL_\tau$ variables for groves is replaced by variables
that are indexed by the $n!$ pairings between the odd nodes and the
even nodes (``odd-even pairings'').  If $\tau$ is such an odd-even pairing, then we define
$$ \X_\tau = \X(\tau) = \prod_{\text{$i$ odd}} \X_{i,\tau(i)} $$
and it turns out to be more convenient to work with
$$ \X'_\tau = \X'(\tau) = (-1)^{\text{\#crosses in $\tau$}}\prod_{\text{$i$ odd}} \X_{i,\tau(i)} $$
where a cross in a pairing $\tau$ is a set of four nodes $a<b<c<d$ such
that $a$ and $c$ are paired with one another and $b$ and $d$ are paired with one another.
The ``$\X$ polynomials'' are in fact integer linear combinations of the $\X'_\tau$'s:
$$ \hPr(\sigma) = \sum_{\text{odd-even pairings $\tau$}} \pd_{\sigma,\tau} \X'_\tau .$$

As with grove partitions, we construct a projection matrix, but with
double-dimer pairings the matrix projects a vector space with basis
vectors indexed by odd-even pairings to a vector space whose basis
vectors are indexed by planar pairings.  Recall that for groves the
projection matrix $\p$ has dimensions $C_n \times B_n$, the $n$\th
Catalan number by the $n$\th Bell number; for double-dimers the
projection matrix $\pd$ has dimensions $C_n \times n!$.  When $n=4$
the projection matrix $\pd$ is

\vspace{3pt}
\newcommand{\n}{$-1\fm$}
\centerline{
\begin{tabular}{r@{$\,\,$}c@{$\,\,$}c@{$\,\,$}c@{$\,\,$}c@{$\,\,$}c@{$\,\,$}c@{$\,\,$}c@{$\,\,$}c@{$\,\,$}c@{$\,\,$}c@{$\,\,$}c@{$\,\,$}c@{$\,\,$}c@{$\,\,$}c@{$\,\,$}c@{$\!\!$}c@{$\!\!$}c@{$\!\!$}c@{$\!\!$}c@{$\!\!$}c@{$\!\!$}c@{$\!\!$}c@{$\!\!$}c@{$\!\!$}c}
                                                          &
\rf{\,{}^1_2\!\mid\!{}^3_4\!\mid\!{}^5_6\!\mid\!{}^7_8\,} &
\rf{\,{}^1_2\!\mid\!{}^3_4\!\mid\!{}^5_8\!\mid\!{}^7_6\,} &
\rf{\,{}^1_2\!\mid\!{}^3_6\!\mid\!{}^5_4\!\mid\!{}^7_8\,} &
\rf{\,{}^1_2\!\mid\!{}^3_8\!\mid\!{}^5_4\!\mid\!{}^7_6\,} &
\rf{\,{}^1_2\!\mid\!{}^3_8\!\mid\!{}^5_6\!\mid\!{}^7_4\,} &
\rf{\,{}^1_4\!\mid\!{}^3_2\!\mid\!{}^5_6\!\mid\!{}^7_8\,} &
\rf{\,{}^1_4\!\mid\!{}^3_2\!\mid\!{}^5_8\!\mid\!{}^7_6\,} &
\rf{\,{}^1_6\!\mid\!{}^3_2\!\mid\!{}^5_4\!\mid\!{}^7_8\,} &
\rf{\,{}^1_6\!\mid\!{}^3_4\!\mid\!{}^5_2\!\mid\!{}^7_8\,} &
\rf{\,{}^1_8\!\mid\!{}^3_2\!\mid\!{}^5_4\!\mid\!{}^7_6\,} &
\rf{\,{}^1_8\!\mid\!{}^3_2\!\mid\!{}^5_6\!\mid\!{}^7_4\,} &
\rf{\,{}^1_8\!\mid\!{}^3_4\!\mid\!{}^5_2\!\mid\!{}^7_6\,} &
\rf{\,{}^1_8\!\mid\!{}^3_4\!\mid\!{}^5_6\!\mid\!{}^7_2\,} &
\rf{\,{}^1_8\!\mid\!{}^3_6\!\mid\!{}^5_4\!\mid\!{}^7_2\,} &
\rf{\,{}^1_2\!\mid\!{}^3_6\!\mid\!{}^5_8\!\mid\!{}^7_4\,} &
\rf{\,{}^1_4\!\mid\!{}^3_6\!\mid\!{}^5_2\!\mid\!{}^7_8\,} &
\rf{\,{}^1_4\!\mid\!{}^3_8\!\mid\!{}^5_2\!\mid\!{}^7_6\,} &
\rf{\,{}^1_4\!\mid\!{}^3_8\!\mid\!{}^5_6\!\mid\!{}^7_2\,} &
\rf{\,{}^1_6\!\mid\!{}^3_2\!\mid\!{}^5_8\!\mid\!{}^7_4\,} &
\rf{\,{}^1_6\!\mid\!{}^3_4\!\mid\!{}^5_8\!\mid\!{}^7_2\,} &
\rf{\,{}^1_6\!\mid\!{}^3_8\!\mid\!{}^5_4\!\mid\!{}^7_2\,} &
\rf{\,{}^1_8\!\mid\!{}^3_6\!\mid\!{}^5_2\!\mid\!{}^7_4\,} &
\rf{\,{}^1_4\!\mid\!{}^3_6\!\mid\!{}^5_8\!\mid\!{}^7_2\,} &
\rf{\,{}^1_6\!\mid\!{}^3_8\!\mid\!{}^5_2\!\mid\!{}^7_4\,}
\\
$\,{}^1_2\!\mid\!{}^3_4\!\mid\!{}^5_6\!\mid\!{}^7_8\,$ & 1& 0& 0& 0& 0& 0& 0& 0& 0& 0& 0& 0& 0& 0&\n&\n& 0&\n& 0&\n& 0& 0& $-2\fm$& 1 \\[2pt]
$\,{}^1_2\!\mid\!{}^3_4\!\mid\!{}^5_8\!\mid\!{}^7_6\,$ & 0& 1& 0& 0& 0& 0& 0& 0& 0& 0& 0& 0& 0& 0& 1& 0&\n& 0& 0& 1& 0& 0& 1&\n       \\[2pt]
$\,{}^1_2\!\mid\!{}^3_6\!\mid\!{}^5_4\!\mid\!{}^7_8\,$ & 0& 0& 1& 0& 0& 0& 0& 0& 0& 0& 0& 0& 0& 0& 1& 1& 0& 0& 0& 0&\n& 0& 1&\n       \\[2pt]
$\,{}^1_2\!\mid\!{}^3_8\!\mid\!{}^5_4\!\mid\!{}^7_6\,$ & 0& 0& 0& 1& 0& 0& 0& 0& 0& 0& 0& 0& 0& 0&\n& 0& 1& 0& 0& 0& 1& 0&\n& 1       \\[2pt]
$\,{}^1_2\!\mid\!{}^3_8\!\mid\!{}^5_6\!\mid\!{}^7_4\,$ & 0& 0& 0& 0& 1& 0& 0& 0& 0& 0& 0& 0& 0& 0& 1& 0& 0& 1& 0& 0& 0& 0& 1& 0       \\[2pt]
$\,{}^1_4\!\mid\!{}^3_2\!\mid\!{}^5_6\!\mid\!{}^7_8\,$ & 0& 0& 0& 0& 0& 1& 0& 0& 0& 0& 0& 0& 0& 0& 0& 1& 0& 1&\n& 0& 0& 0& 1&\n       \\[2pt]
$\,{}^1_4\!\mid\!{}^3_2\!\mid\!{}^5_8\!\mid\!{}^7_6\,$ & 0& 0& 0& 0& 0& 0& 1& 0& 0& 0& 0& 0& 0& 0& 0& 0& 1& 0& 1& 0& 0& 0& 0& 1       \\[2pt]
$\,{}^1_6\!\mid\!{}^3_2\!\mid\!{}^5_4\!\mid\!{}^7_8\,$ & 0& 0& 0& 0& 0& 0& 0& 1& 0& 0& 0& 0& 0& 0& 0&\n& 0& 0& 1& 0& 1& 0&\n& 1       \\[2pt]
$\,{}^1_6\!\mid\!{}^3_4\!\mid\!{}^5_2\!\mid\!{}^7_8\,$ & 0& 0& 0& 0& 0& 0& 0& 0& 1& 0& 0& 0& 0& 0& 0& 1& 0& 0& 0& 1& 0& 0& 1& 0       \\[2pt]
$\,{}^1_8\!\mid\!{}^3_2\!\mid\!{}^5_4\!\mid\!{}^7_6\,$ & 0& 0& 0& 0& 0& 0& 0& 0& 0& 1& 0& 0& 0& 0& 0& 0&\n& 0&\n& 0&\n&\n& 1& $-2\fm$ \\[2pt]
$\,{}^1_8\!\mid\!{}^3_2\!\mid\!{}^5_6\!\mid\!{}^7_4\,$ & 0& 0& 0& 0& 0& 0& 0& 0& 0& 0& 1& 0& 0& 0& 0& 0& 0&\n& 1& 0& 0& 1&\n& 1       \\[2pt]
$\,{}^1_8\!\mid\!{}^3_4\!\mid\!{}^5_2\!\mid\!{}^7_6\,$ & 0& 0& 0& 0& 0& 0& 0& 0& 0& 0& 0& 1& 0& 0& 0& 0& 1& 0& 0&\n& 0& 1&\n& 1       \\[2pt]
$\,{}^1_8\!\mid\!{}^3_4\!\mid\!{}^5_6\!\mid\!{}^7_2\,$ & 0& 0& 0& 0& 0& 0& 0& 0& 0& 0& 0& 0& 1& 0& 0& 0& 0& 1& 0& 1& 0&\n& 1&\n       \\[2pt]
$\,{}^1_8\!\mid\!{}^3_6\!\mid\!{}^5_4\!\mid\!{}^7_2\,$ & 0& 0& 0& 0& 0& 0& 0& 0& 0& 0& 0& 0& 0& 1& 0& 0& 0& 0& 0& 0& 1& 1& 0& 1       \\[2pt]
\end{tabular}}
\vspace{3pt}
We call this matrix $\pd$ the \textbf{projection matrix from odd-even pairings to planar pairings}.

For example, the first row tells us that
\begin{align*}
\hPr(\,{}^1_2\!\mid\!{}^3_4\!\mid\!{}^5_6\!\mid\!{}^7_8\,)
=& \fm
   \X'(\,{}^1_2\!\mid\!{}^3_4\!\mid\!{}^5_6\!\mid\!{}^7_8\,)
 - \X'(\,{}^1_2\!\mid\!{}^3_6\!\mid\!{}^5_8\!\mid\!{}^7_4\,)
 - \X'(\,{}^1_4\!\mid\!{}^3_6\!\mid\!{}^5_2\!\mid\!{}^7_8\,)
 - \X'(\,{}^1_4\!\mid\!{}^3_8\!\mid\!{}^5_6\!\mid\!{}^7_2\,) \\
&- \X'(\,{}^1_6\!\mid\!{}^3_4\!\mid\!{}^5_8\!\mid\!{}^7_2\,)
-2 \X'(\,{}^1_4\!\mid\!{}^3_6\!\mid\!{}^5_8\!\mid\!{}^7_2\,)
 + \X'(\,{}^1_6\!\mid\!{}^3_8\!\mid\!{}^5_2\!\mid\!{}^7_4\,) \\
=& \fm
   \X_{1,2}\X_{3,4}\X_{5,6}\X_{7,8}
 + \X_{1,2}\X_{3,6}\X_{5,8}\X_{7,4}
 + \X_{1,4}\X_{3,6}\X_{5,2}\X_{7,8}
 + \X_{1,4}\X_{3,8}\X_{5,6}\X_{7,2} \\
 &+ \X_{1,6}\X_{3,4}\X_{5,8}\X_{7,2}
-2 \X_{1,4}\X_{3,6}\X_{5,8}\X_{7,2}
 + \X_{1,6}\X_{3,8}\X_{5,2}\X_{7,4}.
\end{align*}

As with the matrix $\p$, we may compute the $\tau$\th column of the
matrix $\pd$ via a sequence of simple combinatorial transformations on
the odd-even pairings.  The prototypical example of the transformation rule is
\begin{center}
\begin{psfrags}
\psfrag{+}[cc][Bc][1][0]{$+$}
\psfrag{-}[cc][Bc][1][0]{$-$}
\psfrag{->}[cc][Bc][1][0]{$\to$}
\includegraphics{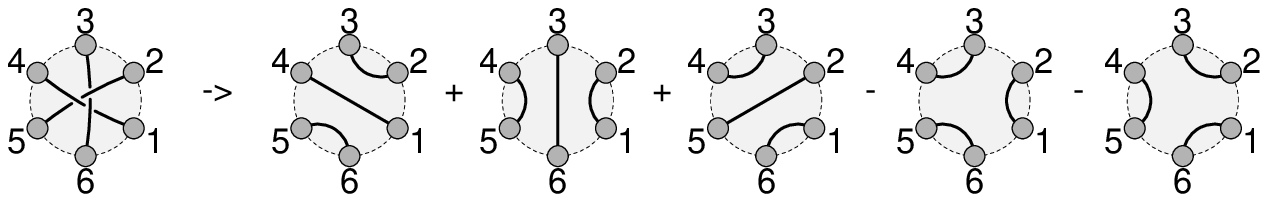}
\end{psfrags}
\end{center}
$$
 \,{}^1_4\!\mid\!{}^3_6\!\mid\!{}^5_2\, \to
 \,{}^1_4\!\mid\!{}^3_2\!\mid\!{}^5_6\,
+\,{}^1_2\!\mid\!{}^3_6\!\mid\!{}^5_4\,
+\,{}^1_6\!\mid\!{}^3_4\!\mid\!{}^5_2\,
-\,{}^1_2\!\mid\!{}^3_4\!\mid\!{}^5_6\,
-\,{}^1_6\!\mid\!{}^3_2\!\mid\!{}^5_4\,
$$
(which is derived in \lref{3Xidentity}).
Compare with the 16\th column, column $\,{}^1_4\!\mid\!{}^3_6\!\mid\!{}^5_2\!\mid\!{}^7_8\,$, in the above matrix.
More generally, suppose that $n\geq 3$ and the odd-even pairing is
$\,{}^{a_1}_{b_1}\!\mid\!{}^{a_2}_{b_2}\!\mid\cdots\mid\!{}^{a_n}_{b_n}\,$,
where the $a_i$'s are odd and the $b_i$'s are even.
Then the transformation rule is
\begin{multline}\label{pair-trans}\tag{Rule~2}
  \,{}^{a_1}_{b_1}\!\mid\!{}^{a_2}_{b_2}\!\mid\!{}^{a_3}_{b_3}\!\mid\rest \to
  \,{}^{a_1}_{b_2}\!\mid\!{}^{a_2}_{b_1}\!\mid\!{}^{a_3}_{b_3}\!\mid\rest
+ \,{}^{a_1}_{b_1}\!\mid\!{}^{a_2}_{b_3}\!\mid\!{}^{a_3}_{b_2}\!\mid\rest
+ \,{}^{a_1}_{b_3}\!\mid\!{}^{a_2}_{b_2}\!\mid\!{}^{a_3}_{b_1}\!\mid\rest\\
- \,{}^{a_1}_{b_2}\!\mid\!{}^{a_2}_{b_3}\!\mid\!{}^{a_3}_{b_1}\!\mid\rest
- \,{}^{a_1}_{b_3}\!\mid\!{}^{a_2}_{b_1}\!\mid\!{}^{a_3}_{b_2}\!\mid\rest
\end{multline}
where in the above ``$\rest$'' represents
${}^{a_4}_{b_4}\!\mid\cdots\mid\!{}^{a_n}_{b_n}\,$ with unchanged
pairings.

For example, when transforming $ \,{}^1_4\!\mid\!{}^3_6\!\mid\!{}^5_8\!\mid\!{}^7_2\,$ we can hold the pair
$\{7,2\}$ fixed:
\begin{align*}
 \,{}^1_4\!\mid\!{}^3_6\!\mid\!{}^5_8\!\mid\!{}^7_2\, \to&
 \,{}^1_6\!\mid\!{}^3_4\!\mid\!{}^5_8\!\mid\!{}^7_2\,
+\,{}^1_4\!\mid\!{}^3_8\!\mid\!{}^5_6\!\mid\!{}^7_2\,
+\,{}^1_8\!\mid\!{}^3_6\!\mid\!{}^5_4\!\mid\!{}^7_2\,
-\,{}^1_6\!\mid\!{}^3_8\!\mid\!{}^5_4\!\mid\!{}^7_2\,
-\,{}^1_8\!\mid\!{}^3_4\!\mid\!{}^5_6\!\mid\!{}^7_2.\\
\intertext{Of these odd-even pairings, the third and fifth ones are planar, but the others require additional applications of Rule~2.
When transforming the first of these terms, if we hold the pair $\{3,4\}$ fixed,
}
 \,{}^1_6\!\mid\!{}^3_4\!\mid\!{}^5_8\!\mid\!{}^7_2\, \to&
 \,{}^1_8\!\mid\!{}^3_4\!\mid\!{}^5_6\!\mid\!{}^7_2\,
+\,{}^1_6\!\mid\!{}^3_4\!\mid\!{}^5_2\!\mid\!{}^7_8\,
+\,{}^1_2\!\mid\!{}^3_4\!\mid\!{}^5_8\!\mid\!{}^7_6\,
-\,{}^1_8\!\mid\!{}^3_4\!\mid\!{}^5_2\!\mid\!{}^7_6\,
-\,{}^1_2\!\mid\!{}^3_4\!\mid\!{}^5_6\!\mid\!{}^7_8\,,
\end{align*}
then the resulting odd-even pairings are all planar.  The other terms above
may be similarly transformed into linear combinations of planar pairings,
and when they are added up, the result is summarized in column $\,{}^1_4\!\mid\!{}^3_6\!\mid\!{}^5_8\!\mid\!{}^7_2\,$ of the projection matrix $\pd$.

\enlargethispage{9pt}
In \sref{ddimers} we prove
\begin{theorem}\label{P2-thm}
  Any odd-even pairing $\tau$ may be transformed into a formal linear
  combination of planar pairings by repeated application of
  \ref{pair-trans}, and the resulting linear combination does not
  depend on the choices made when applying \ref{pair-trans}, so
  that we may write
  $$\tau \to \sum_{\text{\rm planar pairings $\sigma$}} \pd_{\sigma,\tau} \sigma.$$
  For any planar pairing~$\sigma$, these same coefficients $\pd_{\sigma,\tau}$ satisfy the equation
 $$\Pr(\sigma) \frac{\ZWB}{\ZBW} = \sum_{\text{\rm odd-even pairings $\tau$}} \pd_{\sigma,\tau} \X'(\tau)$$
for bipartite circular planar graphs.
\end{theorem}

In \sref{P0P2} we show that the $\pd$ projection matrix of order $n$ is up to signs
embedded in the $\p$ projection matrix of order $2n$:
\begin{theorem} \label{Pt-P2}
 $$ \pd_{\sigma,\tau} = (-1)^{\sigma^{-1} \tau} \p_{\sigma,\tau}$$
where on the left $\sigma$ and $\tau$ denote odd-even pairings, and on the right they
are interpreted as partitions consisting of parts of size $2$, and in the sign
they are interpreted as maps from odd nodes to even nodes, so that $\sigma^{-1}\tau$ is
a permutation on odd nodes, and $(-1)^{\sigma^{-1} \tau}$ is its signature.
\end{theorem}

\subsection{Multichordal SLE connection probabilities}
\label{SLE}

In \sref{SLEsection} we consider the connection probabilities
in the scaling limits of the spanning tree and double-dimer models,
and also of the contour lines in the scaling limit of the discrete
Gaussian free field with certain boundary conditions.

\begin{figure}[b!]
\centerline{\includegraphics[width=0.398\textwidth,angle=90]{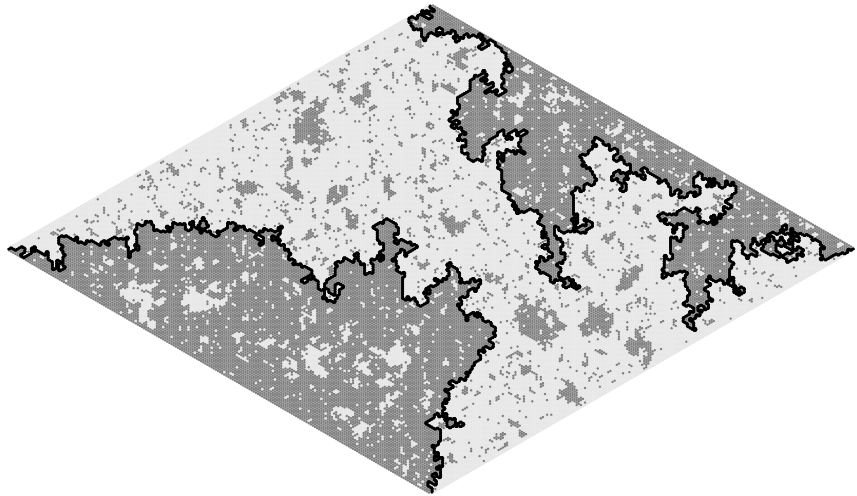}}
\caption{The critical 2D Ising model with mixed spin-up and spin-down
  boundary conditions, where the spins are shown as black and white
  hexagons.  In this configuration there is a connection within the
  white spins connecting the two white boundary segments.  The two
  spin interfaces are shown in bold, and in the scaling limit converge
  to bichordal $\SLE_3$.\label{ising}}
\end{figure}

Connection probabilities of this sort were first studied by Cardy, who
gave (a physics derivation for) an explicit formula for the
probability (in the scaling limit) of a percolation crossing from one
segment of the boundary of a domain to another segment of the boundary
\cite{MR92m:82048}.  Carleson noticed that Cardy's formula is
especially nice when the domain is an equilateral triangle, and this
was one of the insights that led to Smirnov's proof that the percolation
interface converges to $\SLE_6$ \cite{MR1851632}
(see also \cite{camia-newman:proof-converge}).  (The $\SLE_\kappa$
process, which we do not define here, was introduced by Schramm
\cite{MR1776084} and describes the scaling limits of random curves
arising in statistical physics --- see \cite{MR2334202}.)
Arguin and Saint-Aubin
\cite[\S~3]{MR1928682} gave (a physics derivation of) the corresponding
crossing probabilities for spins of the critical Ising model, shown in
\fref{ising}.
As \fref{ising} shows, there are two spin interfaces
that connect the four boundary locations where the boundary conditions
change, there are two ways in which the interfaces may pair up these boundary
points, and the pairing probabilities are computed in \cite[\S~3]{MR1928682}.
Each spin interface,
conditional on the location of the other interface, is distributed
according to a chordal $\SLE_3$ curve \cite{smirnov:ising} within the
domain cut by the other interface, so Arguin and Saint-Aubin
effectively studied a bichordal version of $\SLE_3$.  Kirchhoff's
formula $\Pr(1|2)/\Pr(12)=R_{1,2}$ can, in view of \cite{\LSW}, be
interpreted as giving the connection probabilities in a bichordal
version of $\SLE_8$ (see \hyperref[m-SLE_8]{\S~\ref{m-SLE_8}}).
The corresponding connection probabilities for
other values of $\kappa$ were given in \cite[\S~8.2
and~8.3]{MR2187598} (see also \cite[\S~4.1]{MR2253875}).

In a related vein, there has been much study of the physics of
multiple-strand networks of self-avoiding walks (polymers),
loop-erased random walks, and other random paths that are either known
or believed to be related to $\SLE_\kappa$.  For example, the
``watermelon exponents'' describe the scaling behavior of the
partition function when multiple strands meet at a point in the
interior or on the boundary of a domain
\cite{saleur:saw,duplantier-saleur:dense-saw,duplantier:lerw,duplantier:On}.
The boundary watermelon exponents for $L$ strands can be extracted
from the limiting behavior of the connection probabilities of $L$
strands of (multichordal) $\SLE$ connecting $2L$ boundary points when
$L$ of the boundary points are close to each other and the other $L$
boundary points are close to each other, though of course these
connection probabilities contain more information than is summarized
in the exponents.  This approach does not yield the bulk watermelon
exponents, since to obtain the bulk exponents one would need the
endpoints of the strands to be in the interior of the domain.  The
scaling exponents associated with more complicated multiple-strand
networks (for $2\leq \kappa\leq 8$) have also been derived
\cite{duplantier:On,duplantier:networks,ohno-binder:networks}.  See
\cite{duplantier:conformal-geometry} for up-to-date lecture notes on
the physics of networks of polymers and other types of strands and
their relation to SLE.

\enlargethispage{12pt}
It is only natural to consider the connection probabilities of
multichordal $\SLE_\kappa$, which is
defined in \cite{math.PR/0605159}, and has the property that,
conditional on one curve, the remaining curves are distributed
according to multichordal $\SLE_\kappa$ with fewer curves within the
domain cut by the conditioned curve.  Kozdron and Lawler
\cite{math.PR/0605159} studied multichordal $\SLE_\kappa$ for a fixed
connection topology (when $\kappa\leq 4$).  Cardy exhibited several
discrete models with multiple interfaces that arise naturally in the
context of conformal field theory, for which the interfaces have the
same joint distribution (whose scaling limit is thought to be multichordal SLE)
conditional on a particular connection topology, but for which the
connection topologies have different probabilities \cite{cardy:ade}.
(Thus to discuss connection topology probabilities of multichordal
SLE, one must specify a discrete model whose scaling limit is being
taken.)  Comparing two such models, the connection topology
probabilities in one model have different weights compared to another
such model, but these weights do not depend on the domain or on the
location of the boundary points where the boundary conditions change,
so that the set of all connection topology probabilities in one such
model determines the connection topology probabilities in another such
model.  Dub\'edat \cite{MR2253875} analyzed the general multichordal
$\SLE_\kappa$ connection probabilities (and included special
discussion of the cases $\kappa=2,6,8$), but there still remain
open problems regarding these connection probabilities when the number
of curves is large.

The scaling limit calculations in \sref{SLEsection} yield
explicit formulas for the multichordal $\SLE_\kappa$ connection
probabilities in the cases $\kappa=2,4,8$ (for any number of curves,
for any connection topology).
There is a scaling limit relation \cite{\LSW} between branches of
uniform spanning trees on periodic planar graphs and $\SLE_2$.
Essentially, the scaling limit, as the lattice spacing tends to zero,
of a branch of the uniform spanning tree on a bounded domain, tends to
a random simple curve which is equal in law to $\SLE_2$.  Similarly,
the curve which winds between the uniform spanning tree and its dual
spanning tree (\fref{multi-peano}) converges in the scaling
limit to an $\SLE_8$ \cite{\LSW}.  The double-dimer paths are thought
to have a scaling limit that is given by $\SLE_4$, but this has not
been proved.  However, the contour lines in the discrete Gaussian free
field with certain boundary conditions have been proved to converge to
$\SLE_4$ \cite{math.PR/0605337} in the scaling limit, and using the results of
\cite{math.PR/0605337} we prove in \tref{GFF-pairing} that the
probability distribution of the pairings of these contour lines
coincides with the pairing distribution for the double-dimer model.

\section{Grove partitions}
\label{groves}

\subsection{The meander matrix} \label{meander}

\begin{figure}[htbp]
  \centerline{\includegraphics[width=0.9\textwidth]{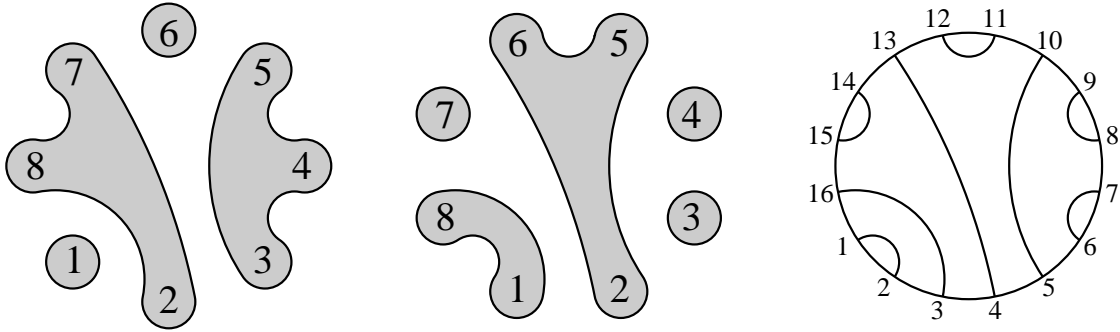}}
\caption{\label{chorddiag} The planar partition $1|278|345|6$, its dual planar partition $18|256|3|4|7$, and their planar chord diagram $\,{}^1_2\!\mid\!{}^{\,\,3}_{16}\!\mid\!{}^{\,\,5}_{10}\!\mid\!{}^7_6\!\mid\!{}^9_8\!\mid\!{}^{11}_{12}\!\mid\!{}^{13}_{\,\,4}\!\mid\!{}^{15}_{14}\,$.}
\end{figure}

Associated to a planar partition $\sigma$ on $n$ nodes
is a \textbf{planar chord diagram}
consisting of $n$ disjoint chords,
winding between the components and the dual components of $\sigma$,
see \fref{chorddiag} for an example.
The chords have in total $2n$ endpoints, one between each node and
adjacent dual node.
A planar chord diagram determines $\sigma$ and vice versa.

There is a natural bilinear form on the space of
formal linear combinations of chord diagrams
that is defined as follows.  Given two chord diagrams $C_1$ and $C_2$,
we can draw them on a sphere, one in the upper
hemisphere and one in the lower hemisphere, so their common boundary
consists of $2 n$ points around the equator.  In the resulting figure
the chords join up to form simple closed loops, each crossing the
equator (at least twice) and such that the set of all loops crosses
the equator a total of $2 n$ times.  Such an object is a \textbf{meander}
of order $n$.  The bilinear form $\langle C_1,C_2\rangle_q$ is defined by
\begin{equation} \label{meander-def}
\langle C_1,C_2\rangle_q = q^{\text{\# loops in meander formed from $C_1$ and $C_2$}}.
\end{equation}
This definition can then be extended linearly to formal linear combinations
of chord diagrams.

The matrix $\M_q$ which has rows and columns indexed by chord diagrams
and has entries $(\M_q)_{C_1,C_2} = \langle C_1,C_2\rangle_q$ is known
as the Gram matrix of the Temperley-Lieb algebra (also called the
meander matrix).  Jones \cite{MR696688} determined the values of $q$
for which $\M_q$ is positive definite and positive semidefinite, Ko
and Smolinsky \cite{MR1105701} determined when $\M_q$ is nonsingular,
and Di Francesco, Golinelli, and Guitter \cite{\DGG} gave an explicit
diagonalization of $\M_q$.  For our results on groves we shall use the
fact that $\lim_{q\to0} \M_q/q$ is nonsingular
\cite[Eqn~(5.18)]{\DGG}, and for our results on the double-dimer model
we shall use the fact that $\M_2$ is nonsingular \cite[Eqn~(5.6)]{\DGG}.

\subsection{The \texorpdfstring{``$L$''}{L} polynomials}

In this subsection we show that $\pu{\sigma}$ is a polynomial in the
$L_{i,j}$ for each planar partition~$\sigma$.

For a planar partition $\sigma$ let $C_\sigma$ be its associated chord
diagram.  We define a bilinear form $\langle,\!\rangle_\g$
on the vector space generated by planar partitions by
\begin{equation}
 \langle\sigma,\tau\rangle_\g =
 \lim_{q\to 0} \frac{\langle C_\sigma,C_\tau\rangle_q}{q} =
 \begin{cases}1, & \text{if meander formed from $C_\sigma$ and $C_\tau$ has one loop,} \\ 0, & \text{otherwise.}\end{cases}
\end{equation}
We let $\M_\g=\lim_{q\to0} \M_q/q$ denote the Gram matrix of
$\langle,\!\rangle_\g$, which is nonsingular.

For a partition~$\sigma$ on $\{1,\dots,n\}$ and graph~$\G$ with nodes
$\{1,\dots,n\}$, let $\G_\sigma$ be the graph $\G$ with nodes
identified according to~$\sigma$, that is $\G_\sigma = \G/\sim$ where
$v\sim v'$ if and only if $v$ and $v'$ are in the same component of~$\sigma$.

If $\sigma$ and $\tau$ are planar partitions and $\G$ is a circular planar graph,
it is not hard to see that
\begin{equation} \label{form-t}
\langle\sigma,\tau\rangle_\g=
 \begin{cases}1, & \text{if a grove of type
$\tau$ in $\G$ gives a spanning tree in $\G_\sigma$}\\[-4pt]&\text{when the vertices
are identified according to $\sigma$,} \\ 0, & \text{otherwise.}\end{cases}
\end{equation}
Indeed, given a grove of $\G$ of type $\tau$, and a ``grove'' of type~$\sigma$
living within the outer face of~$\G$, if their union is a spanning tree,
the path winding around the outside of this spanning tree is a meander.
We can contract to points the components of the ``grove'' of type~$\sigma$,
without changing the topology of the meander.
Another way to say this is that $\langle\sigma,\tau\rangle_\g=1$ if and
only if the number of parts of $\sigma$ and $\tau$ add up to $n+1$ and
the transitive closure of the relation defined by the union of
$\sigma$ and $\tau$ has a single equivalence class.

Let $Z(\tau)$ be the weighted sum of groves of type $\tau$ on $\G$
and $Z_{\G_\sigma}(\tau)$ the corresponding weighted sum on $\G_\sigma$.  In particular
$Z_{\G_\sigma}(\tree)$ is the weighted sum of spanning trees of $\G_\sigma$.
\begin{lemma}\label{ZsigmaZtau}
  For a planar partition~$\sigma$ and circular planar graph~$\G$,
  $$Z_{\G_\sigma}(\tree) = \sum_{\tau}\langle\sigma,\tau\rangle_\g Z(\tau),$$
  where the sum is over all planar partitions.
\end{lemma}

\begin{proof}
  If $\langle\sigma,\tau\rangle_\g=1$ then every grove of $\G$ of type $\tau$, when we
  identify vertices in $\sigma$, becomes a spanning tree of $\G_\sigma$.
  Conversely, any spanning tree of $\G_\sigma$ arises from a unique grove
  of $\G$, and this grove has partition $\tau$ satisfying
  $\langle\sigma,\tau\rangle_\g=1$.
\end{proof}

Note also that $Z_{\G_\sigma}(\unc)=Z(\unc)$, since an uncrossing of $\G$ is
the same as an uncrossing of $\G_\sigma$.

\begin{lemma} \label{glue-L}
  For any graph $\G$ (not necessarily planar) and partition~$\sigma$ (not necessarily planar),
  $Z_{\G_\sigma}(\tree)/Z(\unc)$ is a polynomial in the $L_{i,j}$.
\end{lemma}

\begin{proof}
  Let $S$ be the matrix whose rows index the equivalence classes of
  the relation $\sigma$ and whose columns index the nodes $1,\dots,n$,
  and whose $i,j$-entry is $S_{i,j}=1$ if node~$j$ is in class~$i$,
  and $s_{i,j}=0$ otherwise.

  The response matrix $\Lambda(\G_\sigma)$ for $\G_\sigma$ (see \aref{response})
  is obtained from the response matrix~$\Lambda$ for $\G$ simply as
  \begin{equation}\label{SLS}
  \Lambda(\G_\sigma) = S \Lambda S^T.
  \end{equation}
  That is, putting potential $v_i$ on vertex~$i$ of $\G_\sigma$ is the
  same as putting potential $v_i$ on each node of $\G$ in the
  equivalence class $i$.  The resulting current out of node~$i$ of
  $\G_\sigma$ is the sum of the currents out of nodes of $\G$ in the
  $i$\th equivalence class.
  In particular, entries in $\Lambda(\G_\sigma)$ are sums of entries in~$\Lambda$.  From
  \lref{putree} we have
  $$\det \tilde \Lambda(\G_\sigma)=Z_{\G_\sigma}(\tree)/Z_{\G_\sigma}(\unc),$$
  where $\tilde \Lambda(\G_\sigma)$ is $\Lambda(\G_\sigma)$ with one row and column removed.
  The left-hand side is a polynomial in the $L_{i,j}$, by \eqref{SLS}, and since
  $Z_{\G_\sigma}(\unc)=Z(\unc)$, this concludes the proof.
\end{proof}

Using \lref{ZsigmaZtau}, we can invert matrix $\M_\g$ (since it is nonsingular)
to write $Z(\tau)$ as a rational linear combination of the $Z_{\G_\sigma}(\tree)$ as $\sigma$ varies.
Dividing both sides by $Z(\unc)$ we see that $\pu{\tau}$ is a polynomial
in the $L_{i,j}$ with rational coefficients.  In the next subsection we
prove that (even though $\M_\g^{-1}$ has noninteger entries) the
coefficients of these polynomials are actually integers.

\subsection{Integrality of the coefficients} \label{grove-integer}

We can extend the bilinear form $\langle,\!\rangle_\g$ to work on any
pair of partitions (not necessarily planar), simply by taking the
characterization in \eref{form-t} and dropping the requirement
that the partitions be planar.  Define an ``extended meander matrix''
$\E_\g$ with rows indexed by planar partitions and columns indexed by
all partitions as above, so that
$(\E_\g)_{\sigma,\tau}=\langle\sigma,\tau\!\rangle_\g$.  The matrix
$\E_\g$ has dimensions $C_n \times B_n$, where $C_n$ is the $n\th$
Catalan number, and $B_n$ is the $n$\th Bell number (the number of
partitions on $n$ items).  The matrix $\M_\g$ is the submatrix of
$\E_\g$ containing only the columns for planar partitions.

Note that with this extended definition of $\langle,\!\rangle_\g$,
\lref{ZsigmaZtau} holds for general, nonplanar graphs, provided we
sum over all (not necessarily planar) partitions.  We let
$\vec G$ be the column vector of ``glue variables'', whose entries
are $\vec G_\sigma = Z_{\G_\sigma}(\tree)$ for $\sigma$ a planar partition of $1,\dots, n$.
Let $\vec Z$ be the column vector of partition variables, whose entries
are $\vec Z_\tau = Z(\tau)$ where $\tau$ runs over all partitions.
The extension of \lref{ZsigmaZtau} to the nonplanar setting gives
\begin{equation} \label{glue}
\vec G = \E_\g \vec Z.
\end{equation}

For any not necessarily planar graph $\G$ on $n$ nodes there is an
electrically equivalent complete graph $\K$ on $n$ vertices (in which
every vertex is a node):
it is the graph whose edge $\{i,j\}$ has conductance $L_{i,j}(\G,\No)$.
(The graphs $\G$ and $\K$ are electrically equivalent in the sense that,
when the same voltages are applied to the nodes of $\G$ and $\K$, the
current responses will be the same, i.e., they have the same Dirichlet-to-Neumann matrix.)
Note that each $Z_{\K}(\tau)$ is trivially a polynomial in the $L_{i,j}$,
since each crossing of $\K$ has weight which is a monomial in the $L_{i,j}$.
In fact, we can write this polynomial explicitly: for each part $\lambda$ of $\tau$,
we count the weighted sum of spanning trees of the complete graph
on the vertices in~$\lambda$, and then take the product over the different parts
of $\tau$.  Recalling our definition of $L_\tau$ in \eref{Ltau}, we have that $Z_\K(\tau) = L_\tau$.
For example $$Z_{\K}(1|23|456) = 1\cdot L_{2,3}\cdot(L_{4,5}L_{5,6}+L_{4,5}L_{4,6}+
L_{4,6}L_{5,6}) = \KL_{1|23|456}.$$
We let
$\vec \KL$ be the column vector of partition variables for $\K$:
$\vec \KL_\tau=Z_{\K}(\tau) = \KL_\tau$.

By the preceding equation, for the graph $\K$, \eref{glue} specializes to $\vec G = \E_\g \vec \KL$.
By \lref{glue-L}, for any graph, $\vec G$ is determined by the $L$'s,
so we have $$\vec G  = \E_\g \vec Z =\E_\g \vec \KL$$
for any graph.
For planar graphs, the entries of $\vec Z$ corresponding to nonplanar
partitions are $0$, so $\E_\g \vec Z = \M_\g \vec Z$, and hence
$\M_\g \vec Z = \E_\g \vec \KL$.
Since $\M_\g$ is invertible,
$$\vec Z = \M_\g^{-1} \E_\g \vec \KL.$$
We define $$\p=\M_\g^{-1}\E_\g.$$
We shall see how to compute $\p$ directly, i.e.\ without inverting
$\M_\g$.  The direct computation of $\p$ involves only integer
operations, from which it will follow that the ``$L$ polynomials''
have integer coefficients.

What is the matrix $\p$?  For each planar partition $\sigma$,
the row $\sigma$ tells us $Z_\sigma = \sum_\tau \p_{\sigma,\tau} L_\tau$.  For each
partition $\tau$, the $\tau$\th column gives us a linear combination
of planar partitions $\sum_\sigma \p_{\sigma,\tau} \sigma$ that is equivalent
to $\tau$ in the sense that for any \textit{planar\/} partition $\rho$,
 $$\left\langle\rho,\sum_\sigma \p_{\sigma,\tau} \sigma\right\rangle_\g =
   \sum_\sigma \p_{\sigma,\tau} (\M_\g)_{\rho,\sigma} = (\M_\g \p)_{\rho,\tau} = (\E_\g)_{\rho,\tau} = \langle \rho,\tau\rangle_\g .$$
(We shall see that $\left\langle\rho,\sum_\sigma \p_{\sigma,\tau} \sigma\right\rangle_\g = \langle \rho,\tau\rangle_\g$ for nonplanar partitions $\rho$ too.)


Let us say that two linear combinations of partitions on $n$ items
$\sum_\tau \alpha_\tau \tau$ and $\sum_\tau \beta_\tau \tau$ are
equivalent ($\eqz$) if for any (possibly nonplanar) partition $\rho$
on $n$ items $\sum_\tau \alpha_\tau \langle\rho,\tau\rangle_\g =
\sum_\tau \beta_\tau \langle\rho,\tau\rangle_\g$.  For example,

\begin{lemma} \label{cross}
  $1|234 + 2|134 + 3|124 + 4|123 \eqz 12|34 + 13|24 + 14|23$
\end{lemma}
\begin{proof}
  For any partition $\rho$ with three parts, $\langle
  \rho,\LHS\rangle_\g = 2 = \langle \rho,\RHS\rangle_\g$:
  By symmetry considerations, we need only consider one such parition,
  say $12|3|4$, and $$\langle
  12|3|4,\LHS\rangle_\g = 1+1+0+0 = 0+1+1 = \langle 12|3|4,\RHS\rangle_\g.$$
  For
  partitions $\rho$ with other numbers of parts, $\langle
  \rho,\LHS\rangle_\g = 0 = \langle \rho,\RHS\rangle_\g$,
  since from \eref{form-t} $\langle\rho,\tau\rangle_\g=0$
  whenever $\#\{\text{parts of $\rho$}\} + \#\{\text{parts of $\tau$}\} \neq \#\text{items}+1$.
\end{proof}

As we shall see, this lemma, together with the following two lemmas,
which show how to adjoin new parts and new items to the partitions of
the left- and right-hand sides of an equivalence ($\eqz$), will allow
us to write any partition as an equivalent ($\eqz$) sum of planar
partitions.

\begin{lemma} \label{newpart}
  Suppose $n\geq2$, $\tau$ is a partition of $1,\dots,n-1$, and $\tau
  \eqz \sum_\sigma \alpha_\sigma \sigma$.  Then
  $$\tau|n \eqz \sum_\sigma \alpha_\sigma \sigma|n.$$
\end{lemma}
\begin{proof}
  If $\{n\}$ is a part of $\rho$ then $\langle \rho,\tau|n\rangle_\g
  = 0 = \langle \rho,\RHS\rangle_\g$.  Otherwise $\langle
  \rho,\tau|n\rangle_\g = \langle \rho\smallsetminus n,\tau\rangle_\g =
  \sum_\sigma \alpha_\sigma \langle \rho\smallsetminus n,\sigma\rangle_\g =
  \sum_\sigma \alpha_\sigma \langle \rho,\sigma|n\rangle_\g$.
\end{proof}

If $\tau$ is a partition of $1,\dots, n-1$ and $j\in\{1,\dots,n-1\}$,
we can insert a new item $n$ into the part of $\tau$ that contains
item $j$.  We refer to the resulting partition on $1,\dots,n$ as
``$\tau$ with $n$ inserted into $j$'s part.''

\begin{lemma} \label{split}
  Suppose $n\geq2$, $\tau$ is a partition of $1,\dots,n-1$, $j\in\{1,\dots,n-1\}$, and $\tau
  \eqz \sum_\sigma \alpha_\sigma \sigma$.  Then $$\text{\rm [$\tau$ with $n$
    inserted into $j$'s part]} \eqz \sum_\sigma \alpha_\sigma \text{\rm [$\sigma$ with
    $n$ inserted into $j$'s part]}.$$
\end{lemma}
\begin{proof}
  If $j$ and $n$ are in the same part of partition $\rho$, as well as being in the same part
  of partition $\pi$, then by \eqref{form-t} $\langle \rho,\pi\rangle_\g =0$.
  Thus $\langle
  \rho,\text{[$\tau$ with $n$
    inserted into $j$'s part]}\rangle_\g =0=\langle \rho,\RHS\rangle_\g$.  If $j$ and $n$
  are in separate parts of $\rho$, then let $\rho'$ denote the
  partition obtained from $\rho$ by merging the two parts containing
  $j$ and $n$, and then deleting $n$.  By \eqref{form-t} we have
  \begin{align*}
  \langle
  \rho,\text{[$\tau$ with $n$ inserted into $j$'s part]}\rangle_\g = \langle
  \rho',\tau\rangle_\g = \sum_\sigma \alpha_\sigma \langle
  \rho',\sigma\rangle_\g = \langle \rho,\RHS\rangle_\g
  .& \qedhere\end{align*}
\end{proof}

Lemmas~\ref{cross}, \ref{newpart}, and \ref{split} imply that the
left-hand-side and right-hand-side of transformation Rule~1 are
equivalent ($\eqz$).

\begin{theorem} \label{column}
  For any partition $\tau$, there is an equivalent ($\eqz$) integer linear
  combination of planar partitions $\sum_\sigma \alpha_\sigma \sigma$ (i.e.,
  where each $\alpha_\sigma\in\Z$).
\end{theorem}
\begin{proof}
  We prove this by induction on the number of items in the partition.
  The theorem is true for $n\leq 3$ since each such partition $\tau$
  is already planar (and with \lref{cross}, we see that it is
  also true for $n=4$).  Suppose $\tau$ contains more items.  If
  $\{n\}$ is a part of $\tau$, then we may use the induction
  hypothesis together with \lref{newpart} to find the desired
  linear combination of planar partitions.  Otherwise, item $n$ is in
  the same part as some other item $j$ (if there is more than one
  choice of $j$, it does not matter which one we pick).  By the
  induction hypothesis, we may write $\tau\smallsetminus n \eqz$ an
  integer linear combination of planar partitions, and by
  \lref{split}, we may write $\tau$ as an integer linear
  combination of ``almost planar'' partitions, by which we mean
  partitions that would be planar if the item $n$ were deleted from
  them.

  Next we use Lemmas~\ref{cross}, \ref{newpart}, \ref{split} to
  express an almost-planar partition $\mu$ as an equivalent integer
  linear combination of planar partitions.  We shall use induction on
  the number $k$ of parts of $\mu$ that cross the chord from $j$ to
  $n$.  There is nothing to show if $k=0$, and otherwise we consider
  the part $S$ of $\mu$ closest to $j$ that crosses the chord from $j$
  to $n$.  We let $a=\{i\in S:i<j\}$ and $c=\{i\in S:i>j\}$, both of
  which are nonempty, $b=\text{the part of $\mu$ containing $j$}$, and
  $d=\{n\}$.  Let $a_0\in a$ and $c_0\in c$.  From \lref{cross},
  upon relabeling $1\to a_0$, $2\to j$, $3\to c_0$, and $4\to n$, we get
\begin{align*}
  a_0,c_0|j,n &\eqz a_0|j,c_0,n + j|a_0,c_0,n + c_0|a_0,j,n + n|a_0,j,c_0 - a_0,j|c_0,n - a_0,n|j,c_0.\\
\intertext{Using \lref{split} we may insert the rest of $a$ into the parts containing $a_0$, the rest of $c$ into the parts containing $c_0$, and the rest of $b$ into the parts containing $j$, to obtain}
  a\cup c|b\cup d &\eqz a|b\cup c\cup d + b|a\cup c\cup d + c|a\cup b\cup d + d|a\cup b\cup c - a\cup b|c\cup d - a\cup d|b\cup c.
\end{align*}
  By further application of Lemmas~\ref{newpart} and \ref{split}, we
  may adjoin each of the remaining parts of $\mu$ to the
  left-hand-side (thereby obtaining $\mu$) and to each partition on
  the right-hand-side.  Of the resulting partitions on the
  right-hand-side, the fourth one is planar, and the other partitions
  are almost planar with $k-1$ parts crossing the chords from $n$ to
  the rest of $n$'s part.
\end{proof}

For a partition $\tau$, let $\sum_\sigma \alpha_\sigma \sigma$ be any
linear combination of planar partitions equivalent ($\eqz$) to $\tau$.
Now $\sum_\sigma (\p_{\sigma,\tau}-\alpha_\sigma) \sigma$ lies in the
null-space of $\M_\g$, but $\M_\g$ is nonsingular, so
$\p_{\sigma,\tau}=\alpha_\sigma$ for each $\sigma$.  In particular,
the linear combination promised by \tref{column} is unique, and
gives the $\tau$\th column of $\p$.  This linear combination was obtained
by repeated application of Rule~1, which completes the proof of
\tref{P-thm}, which in turn implies the second part of
\tref{L-R-thm} as a corollary.

We remark that the entries of the matrix $\p$ are all $0$ or $\pm1$
for $n\leq 7$ nodes, but that when $n\geq 8$ other integers appear.

\subsection{The \texorpdfstring{``$R$''}{R} polynomials}

In this subsection we prove the first part of \tref{L-R-thm},
i.e., we show that $\pt{\sigma}$ is a certain polynomial in the
$R_{i,j}$'s for each planar partition $\sigma$.  We start with some
elementary characterizations of the $R_{i,j}$ and $L_{i,j}$ variables
in terms of groves.

\begin{figure}[b!]
\centerline{
\newcommand{\sz}{0.8}
\psfrag{1'}[ct][ct][1][0]{$1$}
\psfrag{2'}[ct][ct][1][0]{$2$}
\psfrag{3'}[ct][ct][1][0]{$3$}
\psfrag{4'}[ct][ct][1][0]{$4$}
\psfrag{5'}[ct][ct][1][0]{$5$}
\psfrag{1}[ct][ct][1][0]{$1$}
\psfrag{2}[ct][ct][1][0]{$2$}
\psfrag{3}[ct][ct][1][0]{$3$}
\psfrag{4}[ct][ct][1][0]{$4$}
\psfrag{5}[ct][ct][1][0]{$5$}
\includegraphics{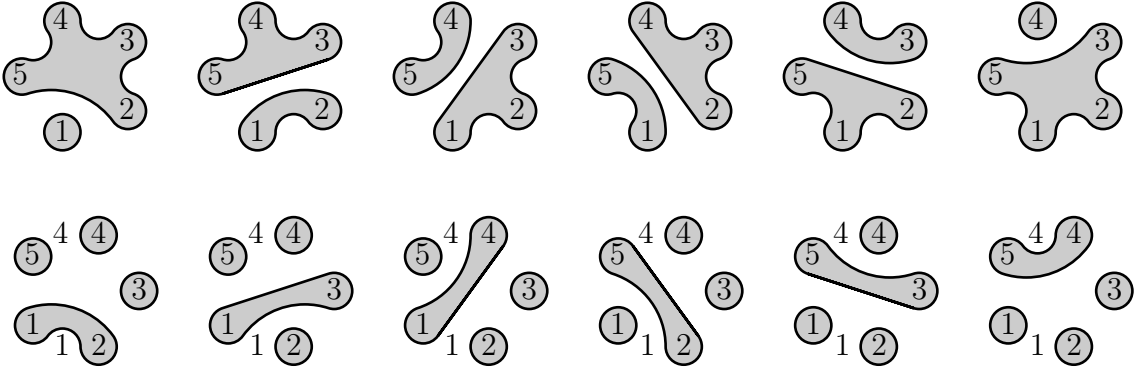}}
\caption{\label{RL} $R_{1,4}$ can be expressed as a sum of 2-part
  partitions in which $1$ and $4$ are in different parts.  For the
  dual resistance $R^*_{1,4}$ in the dual graph, the sum over dual
  partitions becomes a sum over partitions consisting of a doubleton
  part and rest singletons, where the doubleton part separates dual
  nodes $1$ and $4$.  }
\end{figure}

\begin{proposition} \label{R}
  The resistance between $i$ and $j$ is
  $$R_{i,j}=\sum_{A:i\in A, j\in A^c}\pt{A|A^c}.$$
\end{proposition}
\begin{proof}
  If nodes $i$ and $j$ were the only two nodes, then we know from
  Kirchhoff's formula (the case $n=2$) that $R_{i,j}$ is the ratio of the probability that $i$
  and $j$ are in different components of a $2$-component random grove,
  to the probability of a $1$-component grove.  When there are other nodes,
  these $2$-component groves necessarily take the form $A|A^c$ where $i\in A$
  and $j\in A^c$.  (See the top half of \fref{RL}.)
\end{proof}

The following proposition may be derived from Lemma~4.1 of \cite{\CIM},
but for the reader's convenience we provide a short proof.
\begin{proposition} \label{L}
  For $i\neq j$,
  $$L_{i,j} = \pu{i,j|\text{\rm rest singletons}},$$
  i.e., $L_{i,j}=\pu{\sigma}$
  where $\sigma$ is the partition in which every part is a singleton except for the part $\{i,j\}$.
\end{proposition}
\begin{proof}
  Recall from the construction in \tref{column}, that whenever a
  partition $\tau$ is expressed as an equivalent sum of planar
  partitions, each of the partitions has the same number of parts as
  $\tau$.  Since $\sigma$ has $n-1$ parts, and any partition $\tau$
  with $n-1$ parts is already planar, the $\sigma$\th row of $\p$ is
  nonzero only in the $\sigma$\th column, so $\pu\sigma=\vec \KL_\sigma = L_{i,j}$.
\end{proof}

\begin{figure}[t!]
  \centerline{\includegraphics[width=0.35\textwidth]{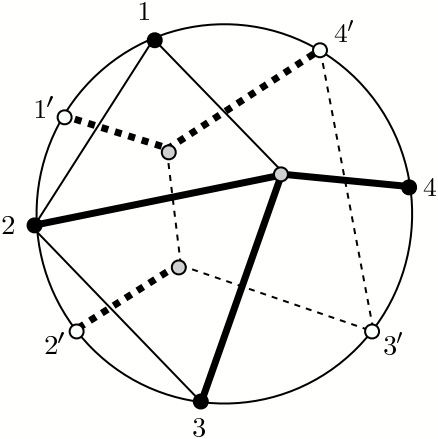}}
\caption{\label{dual} Shown here is a circular planar graph $\G$ with four nodes (shown in black) and one inner vertex, with edges shown as solid lines.  The dual circular planar graph $\G^*$ has four nodes (shown in white) and two inner vertices, with edges shown as dashed lines.  Also shown is a grove of $\G$ of type $1|234$ (the edges of the grove are shown in bold) and its dual grove of $\G^*$ (with edges also shown in bold), which has type $14|2|3$.}
\end{figure}

Any circular planar graph $\G$ has a dual circular planar graph
$\G^*$, as shown in \fref{dual}.  The dual $\G^*$ has an inner
vertex for every bounded face of $\G$, and a node numbered $i$ between
consecutive nodes $i,i+1\bmod n$ of $\G$.  For each edge of $\G$ there
is a dual edge of $\G^*$, whose conductance is, by definition, the
reciprocal of the conductance on the corresponding edge of $\G$.  For
each grove of $\G$ there is a dual grove of $\G^*$ formed from the
duals of the edges of $\G$ not contained in the grove (see
\fref{dual}).  A grove in $\G$ has weight equal to the weight of
its dual grove in $\G^*$, times the product of the conductances of all
edges in $\G$.

\begin{proposition} \label{L*R*}
For the dual graph we have
\begin{equation} \label{L*}
 L^*_{i,j} = \textstyle\frac12 R_{i,j} +  \frac12 R_{i+1,j+1} - \frac12 R_{i,j+1} - \frac12 R_{i+1,j}
\end{equation}
and
\begin{equation} \label{R*}
 R^*_{i,j} =
\doublesum_{\substack{1\leq i'<j'\leq n\\ \text{\rm\ chord $i',j'$ crosses dual chord $i,j$}}}
 L_{i',j'}.
\end{equation}
\end{proposition}

\begin{proof}
  By \pref{R}, $R^*_{i,j}=\sum_{A:i\in A, j\in
    A^c}\Pt^*(A|A^c)$.  But dual nodes $i$ and $j$ are in opposite
  parts of $A|A^c$ if and only if the partition dual to $A|A^c$ is
  $i',j'|\text{rest singleton}$ for some chord $i',j'$ crossing dual
  chord $i,j$ (see \fref{RL}); applying \pref{L} then yields \eqref{R*}.
  Expanding the right-hand-side of~\eqref{L*} using~\eqref{R*} yields
  the left-hand-side of~\eqref{L*}.
\end{proof}

\newpage

\begin{proof}[Proof of first part of \tref{L-R-thm}]
  Since the dual of a partition $\sigma$ is a partition $\sigma^*$ on
  the dual graph, and each grove has the same weight as its dual
  grove, times a constant, we have
  $$\frac{\Pr(\sigma)}{\Pr(\tree)} = \frac{\Pr(\sigma^*)}{\Pr(\text{dual uncrossing})},$$
  which by the second part of \tref{L-R-thm} is an integer-coefficient
  polynomial in the $L^*_{i,j}$'s, which by \pref{L*R*} is an
  integer-coefficient polynomial in the $R_{i,j}/2$'s.
\end{proof}

\section{Double-dimer pairings}\label{ddimers}

\begin{figure}[b!]
  \centerline{\hfill\includegraphics[scale=0.51]{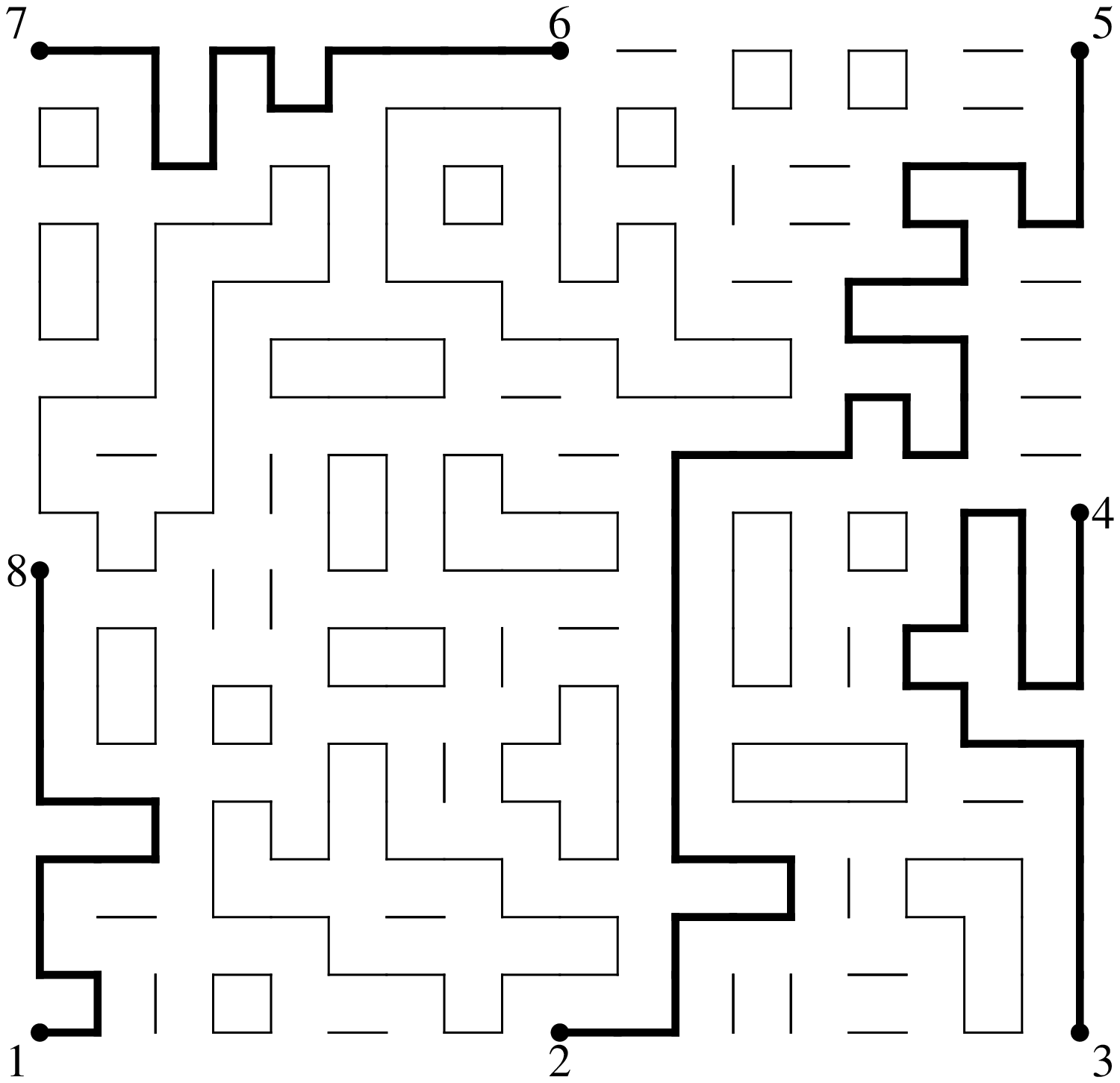}\hfill
\includegraphics[scale=0.51]{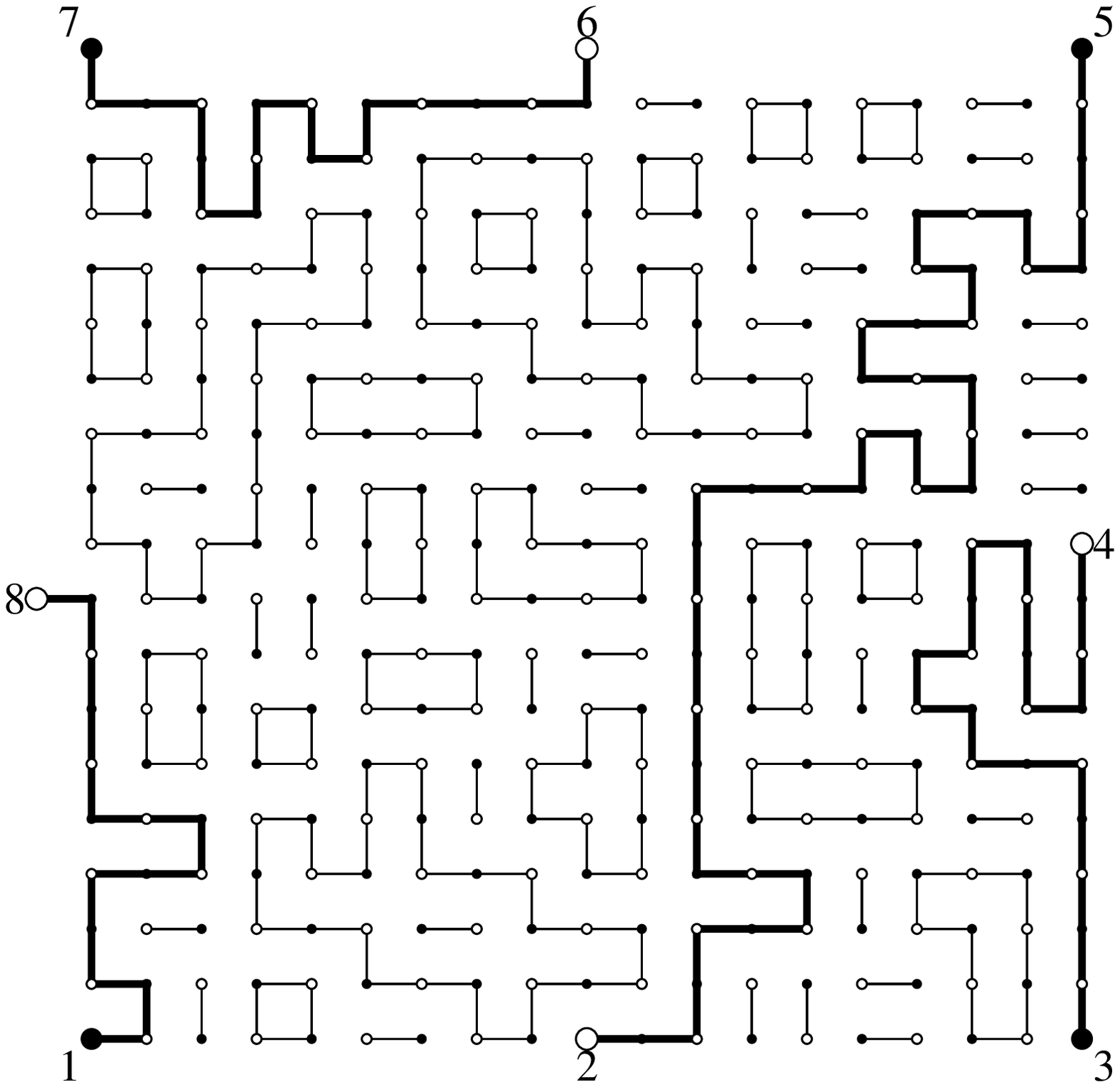}\hfill}
  \caption{\label{ddfig}
    Double-dimer configuration on a rectangular grid with $8$ terminal
    nodes.  In this configuration the pairing of the nodes is
    $\{\{1,8\},\{3,4\},\{5,2\},\{7,6\}\}$, which we write as
    $\,{}^1_8\!\mid\!{}^3_4\!\mid\!{}^5_2\!\mid\!{}^7_6\,$.
    The double-dimer configurations on the graph on the left are in
    one-to-one correspondence with the double-dimer configurations of
    the extended graph on the right, for which the odd terminals are
    colored black and the even terminals are colored white.  }
\end{figure}
Recall that for our double-dimer results we assume the circular planar
graph $\G$ is bipartite, so we may color its vertices black and
white so that each edge connects a black vertex to a white vertex.  It
is convenient to assume that the nodes on the outer face alternate in
color, so that the odd-numbered nodes are black and the even-numbered
nodes are white.  If the graph $\G$ does not satisfy this property, we
can extend $\G$ by adjoining an extra vertex and edge with weight $1$ for each node
that has the wrong color (refer to \fref{ddfig}), and the
double-dimer configurations of the extended graph are in one-to-one
weight-preserving and connection-preserving correspondence with the
double-dimer configurations of the original graph.
Furthermore, the quantities $\ZWB,\ZBW$ as well as the variables
$\X_{i,j}$ are the same on this new graph as they were on
$\G$. Henceforth we assume without loss of generality that the node
colors of $\G$ alternate between black and white.

\subsection{Polynomiality} \label{ddimer-poly}
Since the nodes alternate in color, $\ZBW$ is the weighted sum of dimer covers of $\G$.
Recall that $\ZDD$ is the weighted sum of double-dimer covers of $(\G,\No)$ where
each of the terminal nodes is included in only one of the dimer coverings.
Let $S$ be a balanced subset of nodes, that is, a subset
containing an equal number $k$ of
white and black nodes. Let $\ZD(S)=\ZD(\G\setminus S)$ be the
weighted sum of dimer covers of $\G\setminus S$,
and $\ZD=\ZD(\varnothing)=\ZBW$.
The superposition of a dimer cover of $\G\setminus S$ and a dimer
cover of $\G\setminus S^c$, where $S^c=\No\setminus S$, is
a double-dimer cover of $(\G,\No)$.
In fact, we have
\enlargethispage{12pt}
\begin{lemma}\label{Mlemma}
  $\ZD(S)\ZD(S^c)$ is a sum of double-dimer configurations for all
  connection topologies $\pi$ for which $\pi$ connects no element of
  $S$ to an element of $S^c$.  That is,
\begin{equation}\label{Meqn}
  \ZD(S)\ZD(S^c) =\ZDD\sum_\pi M_{S,\pi} \Pr(\pi),
\end{equation}
  where $M_{S,\pi}$ is $0$ or $1$ according to whether $\pi$ connects
  nodes in $S$ to $S^c$ or not.
\end{lemma}
As a special case, when $S=\varnothing$ we have
$\ZD(\varnothing)\ZD(\No)=\ZDD$, which is closely related to Ciucu's
graph factorization theorem \cite[Thm.~1.2]{ciucu}.  Ciucu showed how
to enumerate dimer coverings in a bipartite graph with bilateral
symmetry.  Dimer coverings of the whole graph correspond to
double-dimer configurations of one half of the graph, with vertices on
the symmetry axis corresponding to nodes, except that there is an
extra factor of $2$ in weight for each path connecting a pair of
nodes, unless the path consists of a single edge.  So the graph
factorization theorem contains a power of~$2$ for each pair of nodes,
and a weight of~$1/2$ for each edge on the symmetry axis, neither of
which appear in \lref{Mlemma}.

Our remaining double-dimer results
only make sense if there are double-dimer configurations of~$\G$,
which implies that there are dimer configurations of~$\G$, so
that~$\G$ has equal numbers of black and white vertices.  After
proving this lemma, we shall henceforth make this assumption.

\begin{proof}[Proof of \lref{Mlemma}]
  Each double-dimer path connecting a pair of nodes separates an even
  number of nodes on its two sides, and since the white and black
  nodes alternate around the boundary, each double-dimer path must go
  from a node to a node of the opposite color.

  Consider the double-dimer cover of $(\G,\No)$ formed by the superposition
  of a dimer cover of $\G\setminus S$ and a dimer cover of $\G\setminus S^c$.
  On the double-dimer path starting from a white node of $S$, every dimer from a
  white vertex to a black vertex is from the second dimer cover and
  every dimer from a black vertex to a white vertex is from the first
  cover.  So such a path necessarily ends at a black vertex in $S$.
  Similarly a path from a black vertex in $S$ necessarily ends at a
  white vertex in $S$.

  Conversely, any double-dimer configuration with no connections from $S$
  to $S^c$ can be decomposed into a dimer cover of $\G\setminus S$ and
  a dimer cover of $\G\setminus S^c$. There are $2^\ell$ choices of
  such a decomposition, $2$ choices for each of the $\ell$ closed loops.
\end{proof}

Before proceeding with the proofs of Theorems~\ref{P2-thm} and~\ref{Z-thm},
we recall some basic facts about
Kasteleyn matrices.  Kasteleyn matrices may be used to enumerate
dimer coverings of any planar graph \cite{\Kasteleyn}, but we review
here just the bipartite case.
A Kasteleyn matrix of an edge-weighted bipartite
planar graph (with a given embedding in the plane) is defined to be a
signed adjacency matrix, $K=(K_{w,b})$ with rows indexed by the white
vertices and columns indexed by the black vertices, satisfying the
following: $K_{w,b}$ is $\pm$ the weight on edge $wb$ and $0$ if there
is no edge.  The signs of the edges are chosen so that around each face there are
an odd number of $-$ signs if the face has $0\bmod 4$ edges and an
even number of $-$ signs if the face has $2\bmod 4$ edges.  Kasteleyn
\cite{\Kasteleyn} showed that every bipartite plane graph with an even
number of vertices has a Kasteleyn matrix, and if there are equal numbers
of black and white vertices then $|\det K|$ is the weighted sum of all
dimer coverings, where the weight of a dimer covering is the product
of its edge weights.

From Kasteleyn theory \cite{\Kasteleyn} \cite{MR1473567}
it is straightforward to compute $\ZD(S)$, though
some work is needed to get all the signs right.  But the signs for
$\ZD(S) \ZD(S^c)$ are simpler, and this product is all we need anyway.
Recall that $\X_{i,j}=\ZD(\{i,j\})/\ZD$.
\begin{lemma}\label{Sdet}
 Let $S$ be a balanced subset of $\{1,\dots,2 n\}$.
 Then
\begin{equation}\label{DS}
\ZD(S)\ZD(S^c)/(\ZD)^2=\det\left[(1_{i,j\in S}+1_{i,j\notin S})\times(-1)^{(|i-j|-1)/2}
\X_{i,j}\right]^{i=1,3\dots,2n-1}_{j=2,4,\dots,2n}.
\end{equation}
\end{lemma}

\begin{proof}
  For convenience we adjoin to the graph $\G$ $2n$ edges along the
  outer face connecting adjacent terminal nodes, and give these edges
  weight $0$ (or weight $\eps$ and then take the limit $\eps\downarrow
  0$).  Given a Kasteleyn matrix of a graph, the signs of edges
  incident to a vertex may be reversed, and each face will still have
  a correct number of minus signs.  For each $i=1,\dots,2n-1$, in
  order, if the edge from node $i$ to node $i+1$ has a minus sign, let
  us reverse the signs of all edges incident to node $i+1$.  Doing
  this ensures that for $1\leq i\leq 2n-1$, the sign of the edge from
  node $i$ to node $i+1$ is positive.  The sign of the remaining
  adjoined edge, from node $2n$ to $1$, will necessarily be $-(-1)^n$
  for the outer face to have a correct number of minus signs.

  Let $(w_1,b_1),\dots,(w_k,b_k)$ be any noncrossing pairing of the
  nodes of $S$, where $w_1,\dots,w_k$ are the white nodes of $S$ and
  $b_1,\dots,b_k$ are the black nodes of $S$.  Let us adjoin edges of
  weight $W$ to the outer face connecting $w_i$ to $b_i$ for $1\leq
  i\leq k$.  To retain the Kasteleyn sign condition, the sign of a new
  edge connecting black node $b$ and white node $w$ will be
  $(-1)^{(|b-w|-1)/2}$.  Let $K_W$ be the Kasteleyn matrix of the
  resulting graph, with the rows and columns ordered so that
  $w_1,\dots,w_k$ are the first $k$ rows and $b_1,\dots,b_k$ are the
  first $k$ columns, and let $K=K_0$ be the corresponding Kasteleyn
  matrix when $W=0$.  Recall that $[x^\alpha]p(x)$ denotes the
  coefficient of $x^\alpha$ in the polynomial $p(x)$.  Then
  $\ZD(S)=\pm [W^k]\det K_W$ and $\ZD=\pm\det K_0$.  But $\det K_W$
  enumerates the weighted matchings of the enlarged graph with all the
  same sign, so
  \begin{equation} \label{zds}
\frac{\ZD(S)}{\ZD} = \frac{[W^k]\det K_W}{[W^0]\det K_W}
 = (-1)^{\sum_{\ell=1}^k \frac{|b_\ell-w_\ell|-1}{2}} \frac{\det K_{\setminus S}}{\det K}
 = (-1)^{\sum_{\ell=1}^k \frac{|b_\ell-w_\ell|-1}{2}} \det [K^{-1}_{b,w}]^{b=b_1,\dots,b_k}_{w=w_1,\dots,w_k}
  \end{equation}
  where $K_{\setminus S}$ denotes the submatrix of $K$ formed from
  deleting the rows and columns from $S$, and the last equality is
  Jacobi's determinant identity.  The special case $S=\{b,w\}$ yields
  $$\X_{b,w}= \ZD(\{b,w\})/\ZD = (-1)^{(|b-w|-1)/2} K^{-1}_{b,w}.$$
  From this equation and~\eqref{zds} we get~\eqref{DS}, except
  possibly for a global sign change.

  Next we compare the sign of $\ZD(S)\ZD(S^c)/(\ZD)^2$ that we get
  from~\eqref{zds} to the sign in~\eqref{DS}.  In the event that
  $S=\{1,\dots,2k\}$, we can take $b_1,\dots,b_k=1,3,\dots,2k-1$ and
  $w_1,\dots,w_k=2,4,\dots,2k$, and for $S^c$ we can take
  $b_1,\dots,b_{n-k}=2k+1,2k+3,\dots,2n-1$ and
  $w_1,\dots,w_{n-k}=2k+2,2k+4,\dots,2n$, so we see that the sign for
  $\ZD(S)\ZD(S^c)/(\ZD)^2$ from~\eqref{zds} agrees with~\eqref{DS}.
  Now suppose that for some $S$ the signs from~\eqref{zds} agree
  with~\eqref{DS}, and we replace one of the nodes $s\in S$ with
  $s+2\notin S$ to get $S'$.  If $s+1\notin S$, the power of $-1$
  in~\eqref{zds} changes by $1$.  If the $s+1\in S$, we may assume $s$
  was paired with $s+1$, and we see that the power of $-1$
  in~\eqref{zds} does not change.  Thus the power of $-1$
  from~\eqref{zds} is opposite for $\ZD(S')\ZD(S'^c)/(\ZD)^2$ compared
  to $\ZD(S)\ZD(S^c)/(\ZD)^2$.  But the product of determinants
  from~\eqref{zds} is $\pm$ the determinant from~\eqref{DS}, and the
  choice of sign is opposite for $(S,S^c)$ and $(S',S'^c)$.
  Thus~\eqref{DS} has the correct global sign for $S'$, and hence by
  induction for any balanced set of nodes.
\end{proof}

By \lref{Sdet} the quantities $\ZD(S)\ZD(S^c)/(\ZD)^2$ are
homogeneous polynomials of degree $n$ (in fact determinants) in the
$\X_{i,j}$, so the matrix $M$ from~\eqref{Meqn} maps the quantities
$\Pr(\pi)\ZDD/(\ZD)^2$ to homogeneous polynomials of degree $n$ in the
$\X_{i,j}$.  We need to show that the matrix $M$ has full rank, that
is, the rank of $M$ is $C_n$, the $n$\th Catalan number, which is the
number of different planar partitions $\pi$.

\begin{lemma}
  $M^T M = \M_2$, the meander matrix $\M_q$ from \eqref{meander-def}
  \cite{\DGG} evaluated at $q=2$.
\end{lemma}

\begin{proof}
  Given two planar matchings $\pi,\sigma$, the integer $\delta_\sigma
  M^T M\delta_\pi$ is the number of subsets $S$ with the property that
  neither $\pi$ nor $\sigma$ connects a node in $S$ to a node in
  $S^c$.

  We draw $\pi$ and $\sigma$ on the sphere, with $\pi$ in the upper
  hemisphere and $\sigma$ in the lower hemisphere, so that their union
  is a meander crossing the equator $2n$ times. Suppose this meander
  has $k$ components.  For each component, the nodes alternate black
  and white, and so the component has the same number of black nodes
  as white nodes. For each component we choose to put all of its nodes
  in $S$ or all of its nodes in $S^c$.  We can construct $2^k$
  possible sets $S$ in this way, and these are exactly the sets $S$
  for which neither $\pi$ nor $\sigma$ connects a node in $S$ to a
  node in $S^c$.

  Thus an entry in $M^T M$ which corresponds to a meander with $k$
  components is $2^k$.
\end{proof}

By \cite[Eqn~(5.6)]{\DGG} the determinant of the meander matrix at $q=2$ is
$$\prod_{i=1}^n(1+i)^{a_{n,i}}$$
where $a_{n,i}$ are certain integers.  In particular it is nonzero,
so $M$ has full rank and we can solve~\eqref{Meqn} for $\Pr(\pi)\ZDD/(\ZD)^2$
in terms of the quantities $\ZD(S)\ZD(S^c)/(\ZD)^2$ which by \lref{Sdet}
are homogeneous polynomials in the $\X_{i,j}$.
This effectively proves \tref{Z-thm}, except for the part
about the coefficients being integers.

\subsection{Integrality of the coefficients} \label{ddimer-integral}

We start by collecting the vectors and matrices we need.
Extend the Gram matrix $M^T M$
when $q=2$ to a $C_n \times n!$ matrix $\E_2$
whose rows are indexed by planar pairings and
whose columns are indexed by all (not necessarily planar)
pairings connecting odd nodes to even nodes.
We have the following matrices and vectors.
\begin{align*}
M =& \text{matrix from balanced subsets $S\subseteq\{1,\dots,2n\}$ to planar pairings $\pi$, from \eqref{Meqn}}\\
\M_2 =& M^T M = \text{the $C_n \times C_n$ meander matrix (with $q=2$) for planar pairings} \\
\E_2 =& \text{$C_n \times n!$ extended meander matrix (with $q=2$) for odd-even pairings}\\
P =& \text{vector of normalized pairing probabilities $P_\pi=\hPr(\pi)=\Pr(\pi) \ZD(\{1,\dots,2n\})/\ZD$}\\
   & \text{indexed by planar pairings $\pi$}\\
D =& MP = \text{vector of products $D_S=\ZD(S)\ZD(S^c)/(\ZD)^2$}\\
   &\phantom{MP=} \text{indexed by balanced subsets $S\subseteq\{1,\dots,2n\}$}\\
\X =& \text{vector of $\X$-monomials indexed by odd-even pairings; $X_\rho = \textstyle\prod_{\{i,j\}\in\rho} X_{i,j}$}\\
\X'=& \text{vector of $\X$-monomials as above with sign $(-1)^{\text{\# crosses of $\rho$}}$}
\end{align*}
Recall that we defined a cross of a pairing $\rho$ to be a set of two
parts $\{a,c\}$ and $\{b,d\}$ of $\rho$ such that $a<b<c<d$.  We
define the parity of an odd-even pairing
$\rho=\,{}^1_{w_1}\!\mid\!{}^3_{w_2}\!\mid\!\cdots\!\mid\!{}^{2n-1}_{w_n}$
to be the parity of the permutation $(w_1/2)(w_2/2)\dots(w_n/2)$.
We shall use the following fact
\begin{lemma}
For odd-even pairings $\rho$,
\begin{equation}\label{crosses}
 (-1)^{\text{\rm parity of $\rho$}}  \prod_{(i,j)\in\rho} (-1)^{(|i-j|-1)/2} = (-1)^{\text{\rm \# crosses of $\rho$}}.
\end{equation}
\end{lemma}
\begin{proof}
  When $\rho = \,{}^1_2\!\mid\!{}^3_4\!\mid\!\cdots$ both sides of the
  above equation are $+1$.  Now suppose we do a transposition,
  swapping the locations of $w$ and $w+2$.  (Such swaps are enough to
  connect the set of odd-even pairings.)  In the event that one of $w$
  or $w+2$ was paired with $w+1$, such a swap will not change the sign
  of the left-hand-side (it changes the sign of exactly one term in
  the product and also the parity of $\rho$), nor does it change the
  number of chords that cross.  Otherwise the sign of the
  left-hand-side does change (exactly two terms of the product change
  sign, and the parity of $\rho$ also changes).  Also, the number of
  chords crossing $w+1$'s chord changes by $0$ or $\pm2$, and the
  chords containing $w$ and $w+2$ now cross if they didn't before, and
  vice versa.
\end{proof}

\begin{lemma}
  $M^T D = \E_2 \X'$
\end{lemma}
\begin{proof}
Recall from \lref{Sdet} that
\begin{align*}
D_S = \frac{\ZD(S)\ZD(S^c)}{(\ZD)^2} &= \det\left[(1_{i,j\in S}+1_{i,j\notin S})\times(-1)^{(|i-j|-1)/2}\X_{i,j}\right]^{i=1,3\dots,2n-1}_{j=2,4,\dots,2n}.\\
\intertext{When we expand the determinant, we get}
D_S &= \sum_{\substack{\text{odd-even pairings $\rho$}\\ \text{$\rho$ does not bridge $S$ to $S^c$}}} \underbrace{(-1)^{\text{parity of $\rho$}}  \prod_{(i,j)\in\rho} (-1)^{(|i-j|-1)/2}
\X_{i,j}}_{\text{\rlap{$=\X'_\rho$ by~\eqref{crosses}}}} \\
\intertext{
Let $\pi$ be a planar pairing.  Upon summing over sets $S$ that are not bridged by $\pi$ we get}
\sum_{\substack{S \subseteq \{1,\dots,2n\}\\ \text{$\pi$ does not bridge $S$ to $S^c$}}} D_S
 &= \sum_{\substack{\text{odd-even pairings $\rho$,}\\ \text{$S$: $\pi,\rho$ do not bridge $S,S^c$}}} \X'_\rho\\ &=
 \sum_{\substack{\text{odd-even}\\\text{ pairings $\rho$}}}
\,\, \sum_{\substack{S:\pi,\rho\text{ do not}\\\text{ bridge $S$}}}\X'_\rho\\
&= \sum_{\substack{\text{odd-even}\\\text{ pairings $\rho$}}} 2^{\text{\# cycles in $\rho\cup\pi$}}\, \X'_\rho.
\end{align*}
The left-hand side is the $\pi$\th row of $M^T D$, and the right-hand side is the $\pi$\th row of $\E_2 \X'$.
\end{proof}

\begin{theorem}
  $ \M_2 P = \E_2 \X'$
\end{theorem}
\begin{proof}
  Since $M P = D$, we have $\M_2 P = M^T M P = M^T D = \E_2 \X'$.
\end{proof}

Since $\M_2$ is invertible, we may define
$$\pd = \M_2^{-1} \E_2.$$
Since $P = \pd \X'$, we can interpret $\pd$ as the matrix
of coefficients of the $\X'$ polynomials: for a given planar pairing
$\sigma$, the $\sigma$\th row of $\pd$ gives the polynomial $\hPr(\sigma)$.
The $\tau$\th column of $\pd$ gives a linear combination of planar pairings
that is in a sense equivalent under $\langle,\!\rangle_2$ to $\tau$: for any
\textit{planar\/} pairing $\rho$,
$$ \left\langle\rho,\sum_\sigma \pd_{\sigma,\tau} \sigma\right\rangle_2 = \sum_\sigma \pd_{\sigma,\tau} (\M_2)_{\rho,\sigma} = (\M_2 \pd)_{\rho,\tau} = (\E_2)_{\rho,\tau} = \langle\rho,\tau\rangle_2.$$
(We shall see $\left\langle\rho,\sum_\sigma \pd_{\sigma,\tau} \sigma\right\rangle_2 = \langle\rho,\tau\rangle_2$ for nonplanar odd-even pairings $\rho$ too.)

We say that two linear combinations of odd-even pairings $\sum_\tau
\alpha_\tau \tau$ and $\sum_\tau \beta_\tau \tau$ are equivalent
($\eqq$) if for any (not necessarily planar) odd-even pairing $\rho$,
we have $\sum_\tau \alpha_\tau \langle\rho,\tau\rangle_q = \sum_\tau
\beta_\tau \langle\rho,\tau\rangle_q$.  Following the approach we used in
\sref{grove-integer} for groves, we will show how to transform an
odd-even pairing $\tau$ into an equivalent ($\eqt$) linear combination
of planar pairings, and thereby compute $\pd$ using only integer operations,
i.e.\ without inverting $\M_2$.

The following key lemma is analogous to \lref{cross}:
\begin{lemma}\label{3Xidentity}
$
 \,{}^1_2\!\mid\!{}^3_4\!\mid\!{}^5_6\,
+\,{}^1_4\!\mid\!{}^3_6\!\mid\!{}^5_2\,
+\,{}^1_6\!\mid\!{}^3_2\!\mid\!{}^5_4\,
\eqt
 \,{}^1_4\!\mid\!{}^3_2\!\mid\!{}^5_6\,
+\,{}^1_2\!\mid\!{}^3_6\!\mid\!{}^5_4\,
+\,{}^1_6\!\mid\!{}^3_4\!\mid\!{}^5_2\,
$
\end{lemma}
\begin{proof}
  For any odd-even pairing $\rho$, we have $\langle\rho,\LHS\rangle_2
  = 12 
  = \langle\rho,\RHS\rangle_2$.
  For example, if $\rho = \,{}^1_4\!\mid\!{}^3_2\!\mid\!{}^5_6\,$ then
  $\langle\rho,\LHS\rangle_2 = 2^2 + 2^2 + 2^2 = 12$, while
  $\langle\rho,\RHS\rangle_2 = 2^3 + 2^1 + 2^1 = 12$.
  For the other odd-even pairings $\rho$ on $\{1,\dots,6\}$,
  it is similarly straightforward to verify $\langle\rho,\LHS\rangle_2
  = \langle\rho,\RHS\rangle_2$ because $2^1+2^1+2^3 = 2^2+2^2+2^2$.
\end{proof}
\lref{3Xidentity} is a special case of \lref{eqq} below.
Though we do not need \lref{eqq}'s extra generality for our
results, some readers may find its proof more instructive.
\begin{lemma}\label{eqq}
If $n$ is a positive integer, then
$$ \sum_{\text{odd-even pairings $\sigma$ on $2n$ items}} (-1)^\sigma \sigma \eqq 0$$
if and only if $q$ is an integer satisfying $0\leq q<n$.
\end{lemma}
\begin{proof}
If $\sigma$ is a permutation on $\{1,2,\dots,n\}$, we can interpret $\sigma$ as the odd-even pairing $\,{}^{\,1}_{2\sigma_1}\!\mid\!{}^{\,3}_{2\sigma_2}\!\mid\!\cdots{}^{\,2n-1}_{2\sigma_n}\,$ and vice versa.
Note that if $\rho$ and $\sigma$ two odd-even pairings, then the cycles in the union of $\rho$ and $\sigma$
are the cycles of $\rho^{-1}\sigma$ when $\rho$ and $\sigma$ are interpreted as permutations, i.e., $\langle\rho,\sigma\rangle_q = q^{\text{\# cycles in $\rho^{-1}\sigma$}}$.
Note also that $(-1)^\sigma = (-1)^n (-1)^\rho (-1)^{\text{\# cycles in $\rho^{-1}\sigma$}}$.
Thus
$$
\left\langle\rho,\sum_{\text{odd-even pairings $\sigma$ on $2n$ items}} (-1)^\sigma \sigma\right\rangle_q
 = (-1)^n(-1)^\rho \sum_{k=1}^n c(n,k) (-q)^k
$$
where the $c(n,k)$ are the unsigned Stirling numbers of the first kind, which count the number of
permutations on $n$ letters that have $k$ cycles.  It is well known (see e.g.\ \cite[Proposition~1.3.4]{MR847717}) that $\sum_k c(n,k) x^k = x(x+1)(x+2)\cdots(x+n-1)$.
\end{proof}

The next lemma is analogous to Lemmas~\ref{newpart} and~\ref{split}:
\begin{lemma} \label{new-pair}
  Suppose $n\geq2$, $\tau$ is an odd-even pairing of $1,2,\dots,2(n-1)$,
  and $\tau \eqq \sum_\sigma \alpha_\sigma \sigma$.  Then
  $$ \tau\!\mid\!{}^{2n-1}_{\,\,2n}\, \eqq \sum_\sigma \alpha_\sigma \sigma\!\mid\!{}^{2n-1}_{\,\,2n}.$$
\end{lemma}
\begin{proof}
  If $\{2n-1,2n\}$ is a part of both odd-even pairings $\rho$ and
  $\pi$, then let $\rho'$ and $\pi'$ denote the odd-even pairings
  obtained by deleting this part.  Then $\langle\rho,\pi\rangle_q = q
  \langle\rho',\pi'\rangle_q$, so $\langle\rho,\LHS\rangle_q =
  \langle\rho,\RHS\rangle_q$.  Now suppose that $\{2n-1,2n\}$ is a
  part of $\pi$ but not $\rho$, and that instead $\{2n-1,w\}$ and
  $\{b,2n\}$ are parts of $\rho$.  Let $\rho'$ denote the odd-even
  pairing obtained from $\rho$ by deleting these two parts and
  replacing them with $\{b,w\}$, and let $\pi'$ be as above.  Then
  $\langle\rho',\pi'\rangle_q = \langle\rho,\pi\rangle_q$, so
  $\langle\rho,\LHS\rangle_q=\langle\rho,\RHS\rangle_q$.
\end{proof}

Lemmas~\ref{3Xidentity} and~\ref{new-pair} imply that the
left-hand-side and right-hand-side of transformation Rule~2 are
equivalent ($\eqt$).

\begin{theorem} \label{pd-column}
  For any odd-even pairing $\tau$, there is an equivalent
  integer-linear combination of planar pairings $\tau\eqt\sum_\sigma
  \alpha_\sigma \sigma$, with $\alpha_\sigma\in\Z$.
\end{theorem}

\begin{proof}
We start with a particular planar pairing,
say $\pi={}^1_2\!\mid\!{}^3_4\!\mid\!{}^5_6\!\mid\!\cdots$.
We have the equivalence $\pi\eqt \pi$.  We shall then make
modifications on the left side, doing adjacent transpositions of even
labels, making the same modification on the right, and then
retransform the right side into a combination of planar pairings.
Eventually we have converted the LHS into $\tau$ and the RHS
is an integer linear combination of planar pairings.

Suppose we swap the labels $b$ and $b+2$.  Assuming that the RHS was a
linear combination of planar pairings, we need to check that after the
swap the RHS can be transformed to again be a linear combination of
planar pairings.  Consider one such planar pairing of the RHS.  If
$b+1$ were paired with one of $b$ or $b+2$, then after the swap the
result is still planar.  Otherwise the three
chords containing $b$, $b+1$, and $b+2$ form three parallel crossings.
These divide the remaining vertices into four contiguous regions, so
that remaining chords only connect vertices in the same region.  After
the $b,b+2$ swap, the three parallel crossings form an asterisk.  The
asterisk can be transformed using Rule~2 (Lemmas~\ref{3Xidentity}
and~\ref{new-pair}).  None of the remaining chords cross the newly
transformed chords, so once they are added back in, the results are
planar.
\end{proof}

For an odd-even pairing $\tau$, let $\sum_\sigma \alpha_\sigma \sigma$ be any
linear combination of planar pairings equivalent ($\eqt$) to $\tau$.
Now $\sum_\sigma (\pd_{\sigma,\tau}-\alpha_\sigma) \sigma$ lies in the
null-space of $\M_2$, but $\M_2$ is nonsingular, so
$\pd_{\sigma,\tau}=\alpha_\sigma$ for each $\sigma$.  In particular,
the linear combination promised by \tref{pd-column} is unique, and
gives the $\tau$\th column of $\pd$.  This linear combination was obtained
by repeated application of Rule~2, which completes the proof of
\tref{P2-thm}, which implies \tref{Z-thm} as a corollary.

\section{Comparing grove and double-dimer polynomials}
\label{P0P2}

As we shall see in this section, the grove partition ``$L$''
polynomials and double-dimer pairing ``$\X$'' polynomials are closely
related to each other.  In fact, the ``$\X$'' polynomials are a
specialization of the ``$L$'' polynomials.
We start with a lemma about groves.

\begin{lemma}
\label{drop-singletons}
  If a planar partition $\sigma$ contains only singleton and doubleton parts,
  and $\sigma'$ is the partition obtained from $\sigma$ by deleting all
  the singleton parts, then the ``$L$'' polynomials for $\sigma$ and
  $\sigma'$ are lexicographically equal.  (We say ``lexicographically
  equal'' rather than ``equal'' because the formulas look the same
  even though the underlying ``$L$ variables'' represent different
  electrical quantities.)
\end{lemma}
For example, $\pu{13|2|4|5|68|7|9} = L_{1,3}L_{6,8} - L_{1,6}L_{3,8}$
and $\pu{13|68} = L_{1,3}L_{6,8} - L_{1,6}L_{3,8}$, but the $L$'s mean different
things in these two equations, since they are the current responses when
different numbers of nodes are held at zero volts.
Note that for general partitions we cannot simply drop the singleton
parts.  For example, $\pu{123} = L_{1,2}L_{2,3} + L_{1,2} L_{1,3} +
L_{1,3} L_{2,3}$, while $\pu{123|4}$ has these three terms together
with a fourth term, namely $L_{1,3} L_{2,4}$.
\begin{proof}
  Recall the algorithm that computes columns of the projection matrix
  by finding for a given partition~$\tau$ an equivalent linear
  combination of planar partitions (the proof of \tref{column}).
  Each partition in the result will
  have the same number of parts as $\tau$, and every singleton of
  $\tau$ is a singleton of each partition in the combination.  If we
  are only interested in the $\sigma$\th row of the projection matrix,
  where $\sigma$ has $k$ parts of size $2$ and rest singletons, then we
  need only consider columns $\tau$ for which $\tau$ has $n-k$ parts.
  If $\tau$ has any part of size $3$ or more, then $\tau$ has more
  singletons that $\sigma$, and it does not contribute to $\sigma$'s row.
  The only partitions $\tau$ that contribute have $k$ parts of size
  $2$, and the same singletons that $\sigma$ does.  For such partitions
  $\tau$, we may drop the singleton parts, find the equivalent linear
  combination of planar partitions, and then re-adjoin the singleton
  parts.  Thus the computation of $\pu{\sigma}$ precisely mirrors the
  computation of $\pu{\sigma'}$.
\end{proof}

The following theorem is a reformulation of \tref{Pt-P2}.

\begin{theorem}
  If a planar partition $\sigma$ contains only doubleton parts, and we make
  the following substitutions to the grove partition polynomial $\pu{\sigma}$:
$$
 L_{i,j} \to \begin{cases} 0, &\text{if $i$ and $j$ have the same parity,} \\
     (-1)^{(|i-j|-1)/2} \X_{i,j}, &\text{otherwise,}\end{cases}
$$
then the result is $(-1)^\sigma$ times
the double-dimer pairing polynomial $\hPr(\sigma)$,
when we interpret $\sigma$ as a pairing, and $(-1)^\sigma$ is the
signature of the permutation $\sigma_1,\sigma_3,\dots,\sigma_{2n-1}$.
\end{theorem}
\begin{proof}
  Consider the computation of $\pu{\sigma}$ when the planar
  partition~$\sigma$ contains only doubleton parts.  When we express a
  partition~$\tau$ as a linear combination of planar partitions, any
  singleton parts of $\tau$ show up in each planar partition (with
  nonzero coefficient), so if~$\tau$ contains singleton parts,
  $\p_{\sigma,\tau}=0$.  When we transform partitions according to
  Rule~1, the number of parts is conserved, so partitions that end up
  contributing to the $\sigma$\th row of the projection matrix will
  contain only doubleton parts.  Rule~1 transforms a partition into a
  combination of six partitions, but the first four of them have the
  wrong part sizes to contribute to $\p_{\sigma,\tau}$, so we may keep
  only keep the last two.  In other words, we may use the abbreviated
  transformation rule
  $$ 13|24 \to -12|34 - 14|23.$$
  Let us transform $14|25|36$ using this abbreviated rule:
  \begin{align*}
    14|25|36 &\to  - 12|45|36 - 15|24|36\\
    &\to  - 12|45|36 + 15|23|46 + 15|26|34 \\
    &\to  - 12|45|36 - 14|23|56 - 16|23|45 - 16|25|34 - 12|56|34.
    \intertext{Compare this to the corresponding rule for
      double-dimer pairings:}
    14|25|36
    &\to  + 12|45|36 + 14|23|56 - 16|23|45 + 16|25|34 - 12|56|34.
  \end{align*}
  The transformation rules are the same except for a sign factor equal
  to the signature of the pairing.  Thus
  $\pd_{\sigma,\tau} = (-1)^{\sigma^{-1} \tau} \p_{\sigma,\tau}$.
  Recall now that the double-dimer
  projection matrix gave coefficients for monomials weighted by
    $(-1)^{\text{\# pairs that cross}}$, and that this sign was
  by~\eqref{crosses} equal to the signature of the pairing times the
  product over pairs $(i,j)$ in the pairing of $(-1)^{(|i-j|-1)/2}$.
\end{proof}

\section{Scaling limits and multichordal SLE}\label{SLEsection}

The results in the preceding sections are all exact results that hold
for any planar graph.  In this section we study the asymptotic
behavior of the partition probabilities when the graphs approximate a
fine lattice restricted to a domain in the plane.  We consider
two different limits for groves; the first one gives
multichordal loop-erased random walk, which in the limit converges to
multichordal $\SLE_2$, while the second one in the limit converges to
multichordal $\SLE_8$.  We also obtain the limiting pairing
probabilities of the paths in the double-dimer model, and the limiting
connection probabilities for the contour lines of the discrete
Gaussian free field with certain boundary conditions, and in
particular show that they are equal.  The latter of these two models
converges to multichordal $\SLE_4$ \cite{math.PR/0605337}
\cite{schramm-sheffield:personal}.
\enlargethispage{12pt}

Dub\'edat \cite{MR2253875} has also studied the pairing
probabilities in multichordal $\SLE_\kappa$.  We only treat the cases
$\kappa=2,4,8$, while Dub\'edat's calculations are relevant for the
continuous range $0\leq\kappa\leq 8$.  On the other hand, our
calculations are more explicit.  For example, Dub\'edat's solution is
actually a multiparameter family of solutions that solve a certain
PDE.  (It is not known \textit{a priori\/} that the solutions in
\cite{MR2253875} span the solution space, though for a given
number of strands
this can be verified \textit{a posteriori}.)
Each of these solutions is relevant to multichordal
$\SLE_\kappa$, but when there are more than two strands, it is not
clear which solution to the PDE is the ``canonical solution'' that
describes the behavior of a given discrete model such as spanning
trees.  Our calculations avoid that issue entirely by working directly
with the discrete models.
\enlargethispage{12pt}

\subsection{Multichordal loop-erased random walk and \texorpdfstring{$\SLE_2$}{SLE(2)}} \label{sle2}

\begin{figure}[b!]
  \centerline{\includegraphics[width=0.49\textwidth]{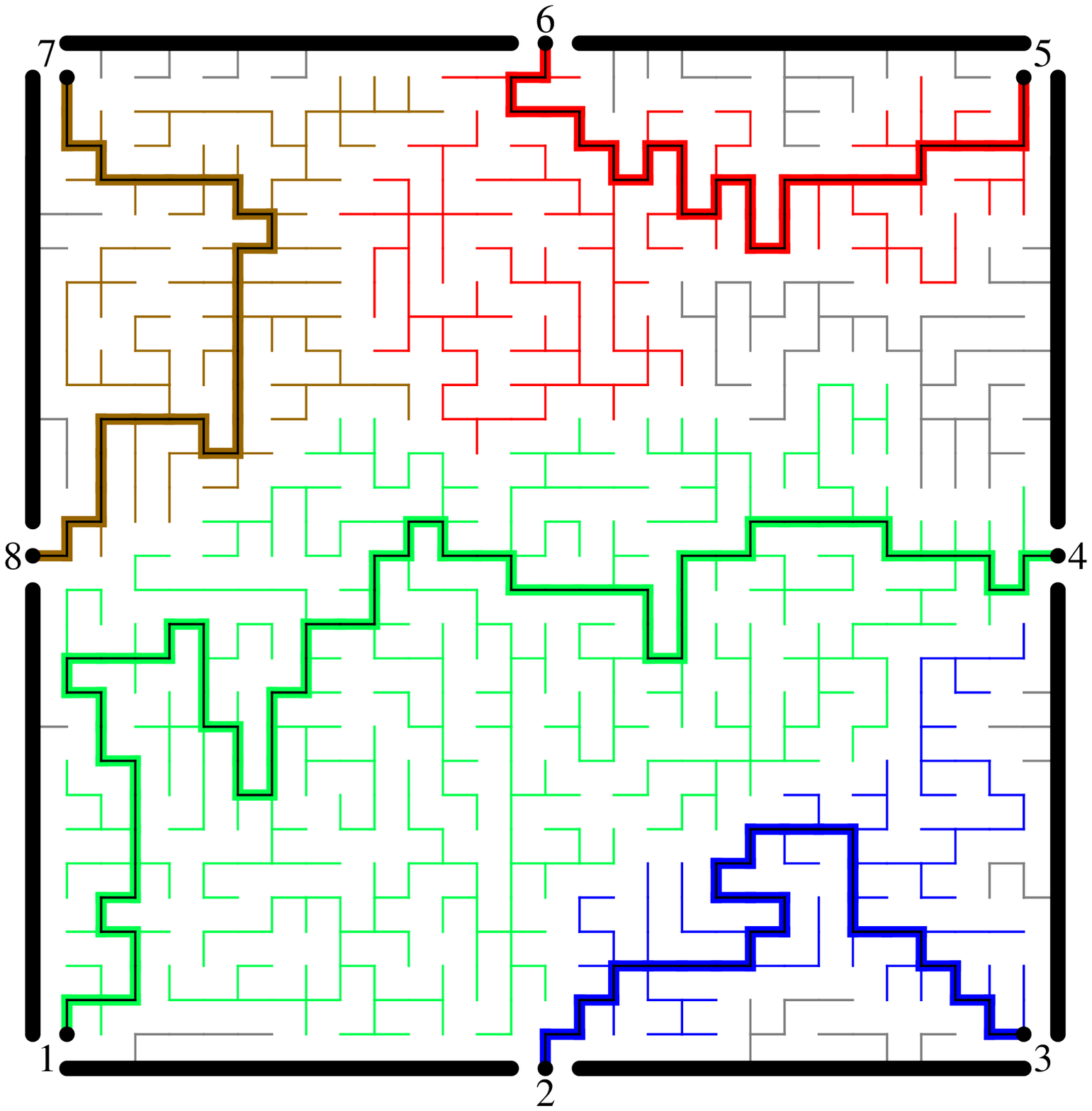}}
\caption{\label{multi-lerw} Multichordal loop-erased random walk: Shown here is a random grove with 16 nodes, conditioned to have the long extended (and unnumbered) nodes in singleton parts,
and the short (numbered) nodes each in a part of size 2.  In this case the pairing is $\,{}^1_4\!\mid\!{}^3_{2}\!\mid\!{}^5_6\!\mid\!{}^7_8\,$,
and in the scaling limit for this geometry this pairing occurs with probability $\frac{48777}{965776}-\frac{1135 \sqrt{2}}{60361} \doteq 0.0239$.
Within the grove the paths connecting the numbered nodes are shown in bold.  Each of these connecting paths is, conditional upon the other paths, the loop-erasure of a random walk started at one endpoint and conditioned to hit the other endpoint before hitting either the boundary or the other paths.
}
\end{figure}

We consider multichordal $\SLE_2$ on a domain $U$,
with $m$ chords connecting $2 m$ points.
To get $\SLE_2$ as a limit of groves, we
consider groves on finer and finer planar grids approximating~$U$,
with $n=4 m$ nodes, where each of the
$2m$ $\SLE$-endpoints is associated to a node, and each of the $2m$ intervals
between the $\SLE$-endpoints is wired together and associated to a node.
(See \fref{multi-lerw}.)
The ``interval'' nodes are required to be
in singleton parts of the corresponding grove partition, and the point nodes are required to be in doubleton parts.  By
\lref{drop-singletons}, to compute the $\SLE_2$ connection
probabilities it is sufficient to compute the current
responses amongst the point nodes (when the interval nodes are
grounded), and substitute these values into the ``$L$'' polynomials
associated with pairings of these $2m$ point nodes.

By conformal invariance, we may without loss of generality fix $U$ to be the upper
half-plane, with point nodes at $x_1,\dots,x_{2m}\in\R$, and we assume
that no point node is at $\infty$.  We approximate the upper half-plane with
the upper half cartesian lattice $\eps\Z\times\eps\N$, and round
$x_i$ to the nearest lattice point $x_i^{(\eps)}\in \eps\Z\times\eps\N$,
and let each edge have unit resistance.
The current response between distinct nodes $x_i^{(\eps)}$ and $x_j^{(\eps)}$
is (assuming the nodes are non-adjacent) the voltage at the vertex one lattice
spacing above $x_j^{(\eps)}$ when $x_i^{(\eps)}$ is at one volt and the remainder
of the real axis is at zero volts, which is the probability that a random walk
started one lattice spacing above $x_j^{(\eps)}$ ends up at
$x_i^{(\eps)}$ when it first hits the real axis, which is
\begin{equation}\label{Lij-lerw}
L_{i,j}=(1+o(1))\frac{1}{\pi} \frac{\eps^2}{(x_i-x_j)^2},
\end{equation}
see e.g.\ \cite[Chapter~III \S 15]{\Spitzer}, where the $o(1)$ term goes to $0$ as $\eps$ goes to $0$.

When we consider the ratios of these probabilities, the $\eps^2/\pi$ factors
drop out and we can then take the limit $\eps\to0$.
We find then that $L_{i,j}$ is
inversely proportional to the square of the distance between points $x_i$ and $x_j$.

In the bichordal case, the normalized probabilities are
\begin{align*}
 \pu{12|34}=&L_{12|34} - L_{13|24} \\
 \pu{14|23}=&L_{14|23} - L_{13|24}.
\end{align*}
In the $\eps\to0$ limit, the unnormalized probability (conditional on there being two chords) is
\begin{align*}
\Pr(14|23)=\frac{\pu{14|23}}{\pu{14|23}+\pu{12|34}}
 &\to \frac{\frac{1}{(x_1-x_4)^2}\frac{1}{(x_2-x_3)^2} - \frac{1}{(x_1-x_3)^2}\frac{1}{(x_2-x_4)^2}}{\frac{1}{(x_1-x_2)^2}\frac{1}{(x_3-x_4)^2} + \frac{1}{(x_1-x_4)^2}\frac{1}{(x_2-x_3)^2}  - 2 \frac{1}{(x_1-x_3)^2}\frac{1}{(x_2-x_4)^2}}. \\
\intertext{If we let $s$ denote the cross ratio $s=(x_4-x_3)(x_2-x_1)/[(x_4-x_2)(x_3-x_1)]$,
this limiting probability can be written as}
 &= \frac{2 s^3-s^4}{1 - 2s + 4 s^3 - 2 s^4},
\end{align*}
which agrees with the known formula for bichordal $\SLE_2$, see \cite[\S~8.2 and~8.3]{MR2187598} \cite{schramm-wilson} \cite[\S~4.1 and~4.2]{MR2253875}.

For trichordal $\SLE_2$, the relevant polynomials are
\begin{align*}
 \pu{12|34|56}=& -L_{15|26|34}+L_{12|56|34}+L_{14|26|35}+L_{15|24|36}\\&-L_{14|25|36}+L_{13|25|46}-L_{12|35|46}-L_{13|24|56} \\
 \pu{16|23|45}=&-L_{16|24|35}+L_{14|26|35}+L_{15|24|36}-L_{14|25|36}\\&+L_{16|23|45}-L_{13|26|45}-L_{15|23|46}+L_{13|25|46} \\
 \pu{14|23|56}=&L_{15|24|36}-L_{14|25|36}-L_{15|23|46}+L_{13|25|46}+L_{14|23|56}-L_{13|24|56} \\
 \pu{25|16|34}=&L_{16|25|34}-L_{15|26|34}-L_{16|24|35}+L_{14|26|35}+L_{15|24|36}-L_{14|25|36} \\
 \pu{36|12|45}=&L_{14|26|35}-L_{12|46|35}-L_{14|25|36}-L_{13|26|45}+L_{12|36|45}+L_{13|25|46}.
\end{align*}

\subsection{Grove Peano curves and multichordal \texorpdfstring{$\SLE_8$}{SLE(8)}} \label{m-SLE_8}

\begin{figure}[b!]
  \centerline{\includegraphics[width=0.48\textwidth]{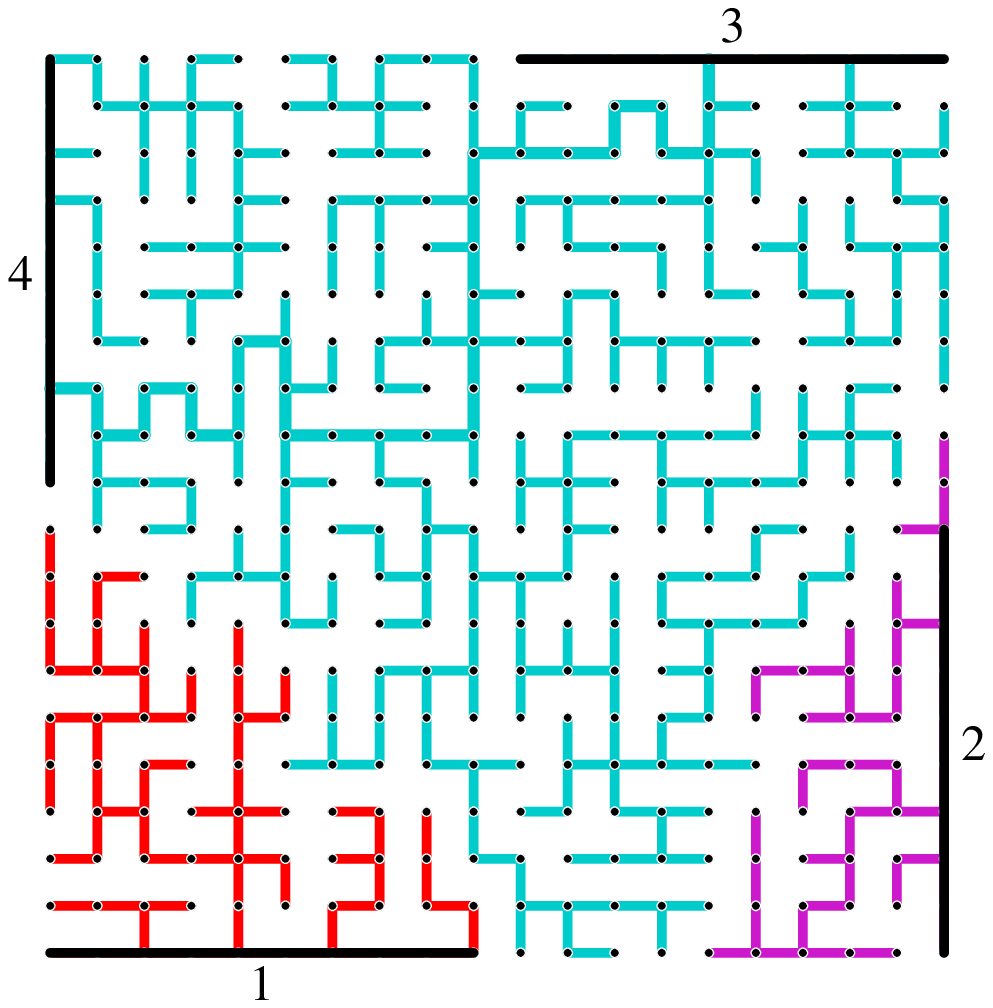}\hfil\includegraphics[width=0.48\textwidth]{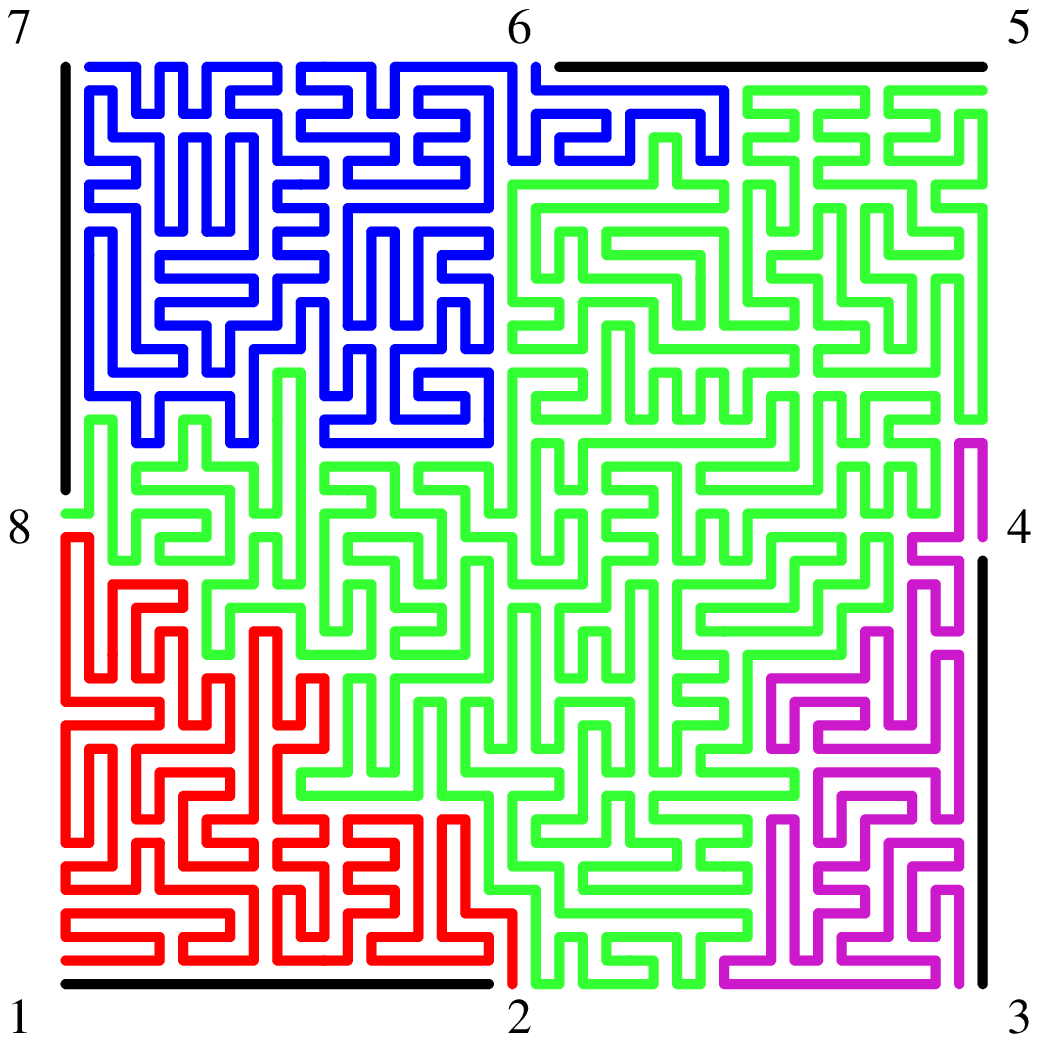}}
\caption{\label{multi-peano}
The multichordal Peano curves (which converge to $\SLE_8$) associated to
a grove with four nodes (the extended black segments).
In this case the grove partition is $1|2|34$, which gives the Peano curve pairing $\,{}^1_2\!\mid\!{}^3_4\!\mid\!{}^5_8\!\mid\!{}^7_6\,$.
In the scaling limit for this geometry these connections occur with probability $(2-\sqrt2)/8 \doteq 0.0732$.
}
\end{figure}

We can do a similar calculation for multichordal $\SLE_8$ containing
$m$ chords connecting $2m$ points.  Lawler, Schramm, and Werner \cite{MR2044671}
showed that on a fine grid in which part of the boundary is wired and part is
free, the path separating a random spanning tree from its dual spanning tree
converges to chordal $\SLE_8$ connecting the two boundary points where
the boundary conditions change from wired to free.  If the boundary conditions
alternate between wired and free and back $m$ times, there will be $m$ paths
separating a random grove from its dual grove which connect the $2m$ points
where the boundary conditions change from wired to free, and it is not difficult
to see that each such path, conditional upon the other paths, converges to $\SLE_8$.
Here the
groves have $n=m$ nodes, which correspond to ``interval'' nodes between
the $(2i-1)$\th and $(2i)$\th endpoints of the $\SLE_8$ strands.
(See \fref{multi-peano}.)
We again consider $U$ to be the upper half-plane and $x_1,\dots,x_{2m}$
the endpoints of the $m$ $\SLE_8$ strands.  Let $I_i$ be the interval
node whose endpoints are $x_{2i-1}$ and $x_{2i}$.  Now $L_{i,j}$ is
the current into $I_j$ for the harmonic function which is $1$ on
$I_i$, zero on the other $I_{j'}$ and has Neumann boundary conditions
elsewhere.

For $m=2$, the event that $x_1,x_4$ are connected and $x_2,x_3$
are connected is the event that the grove has one tree connecting
nodes $1$ (the interval $I_1=[x_1,x_2]$) and $2$ (the interval $I_2=[x_3,x_4]$),
which has probability $\pu{12}/(\pu{12}+\pu{1|2})=L_{1,2}/(L_{1,2}+1)$, where $L_{1,2}$ is the
modulus of the rectangle which is the conformal image of $U$ and whose
vertices are the images of $x_1,x_2,x_3,x_4$.  By the
Schwarz-Christoffel formula, the Riemann map from $U$ to the rectangle
is
$$f(z) = \int_{0}^z \frac{dw}{i (w-x_1)^{1/2}(w-x_2)^{1/2}(w-x_3)^{1/2}(w-x_4)^{1/2}}.$$
Using $f(z)$ we may express
$$L_{1,2} = \frac{f(x_2)-f(x_1)}{-i(f(x_3)-f(x_2))}.$$

For $m=3$, an interesting special case is where the domain is the
regular hexagon and the three nodes are three nonadjacent sides of the
hexagon.  Recall now the dual electric network with its dual grove,
which was illustrated in \fref{dual}.  Here the dual electric
network also converges to the regular hexagon, but with extended nodes
on the other three nonadjacent sides of the hexagon.  Thus the
probability that the primal grove has type $1|2|3$ equals the
probability that the dual grove has type $1|2|3$, but this latter
event is the event that the primal grove has type $123$, so
$\Pr[1|2|3]=\Pr[123]$.  From this we see
$1=\pu{1|2|3}=\pu{123}=L_{1,2}L_{2,3}+L_{1,2}L_{1,3}+L_{1,3}L_{2,3}$,
and hence by symmetry $L_{1,2}=L_{1,3}=L_{2,3}=1/\sqrt{3}$.  We can
now compute the other normalized probabilities
$\pu{12|3}=\pu{13|2}=\pu{23|1}=L_{2,3}=1/\sqrt{3}$.  Thus the
unnormalized partition probabilities are
$\Pr[1|2|3]=\Pr[123]=1/(2\times 1+3\times 1/\sqrt{3})=2-\sqrt3 \doteq 0.268$ and
$\Pr[12|3]=\Pr[13|2]=\Pr[1|23]=2/\sqrt3-1 \doteq 0.154$.  (The
corresponding calculation for the connection probabilities for
percolation ($\SLE_6$) on a regular hexagon was carried out by
Dub\'edat \cite[\S~4.4]{MR2253875}, but while the $\SLE_8$ calculations
here are almost trivial in retrospect, the corresponding $\SLE_6$
calculations are much more intricate.)  The scaling limit of the $L$'s
can similarly be worked out using duality for the example in
\fref{multi-peano}: $L_{i,i+1}=1/2$ and
$L_{i,i+2}=1/\sqrt{2}-1/2$.  The response matrix for the regular
$2n$-gon with sides that are alternately free and wired has a nice
formula that is given in \cite{math.CA/0703826}.

\begin{figure}[b!]
\centerline{
\newcommand{\sz}{0.8}
\psfrag{fa1}[ct][ct][\sz][0]{$f(\!x_5\!)$}
\psfrag{fa2}[lt][lt][\sz][0]{$f(\!x_7\!)$}
\psfrag{fa3}[lt][lt][\sz][0]{$f(\!x_9\!)$}
\psfrag{fa4}[lt][lt][\sz][0]{$f(\!x_1\!)$}
\psfrag{fb1}[rt][rt][\sz][0]{$f(\!x_6\!)$}
\psfrag{fb2}[rt][rt][\sz][0]{$f(\!x_8\!)$}
\psfrag{fb3}[rt][rt][\sz][0]{$f(\!x_{10}\!)$}
\psfrag{fb4}[ct][ct][\sz][0]{$f(\!x_2\!)$}
\psfrag{fc1}[cb][cb][\sz][0]{$f(\alpha_1)$}
\psfrag{fc2}[cb][cb][\sz][0]{$f(\alpha_2)$}
\psfrag{fc3}[cb][cb][\sz][0]{$f(\alpha_3)$}
\psfrag{fa5}[cb][cb][\sz][0]{$f(x_3)$}
\psfrag{fb5}[cb][cb][\sz][0]{$f(x_4)$}
\includegraphics[width=0.5\textwidth]{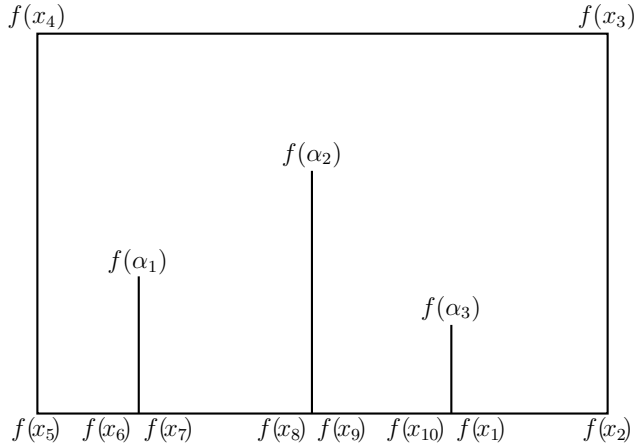}}
\caption{\label{slit-rect}
The function $f(z)=\int_0^z\frac{\prod_{\ell=1}^{m-2}(w-\alpha_\ell)\,dw}{i\prod_{\ell=1}^{2m} (w-x_\ell)^{1/2}}$ maps the upper half-plane to the slit rectangle,
for suitable choices of the $\alpha_\ell$'s.
In the above figure $L_{[x_3,x_4],[x_7,x_8]} = i \frac{f(x_8)-f(x_7)}{f(x_3)-f(x_2)}$.
}
\end{figure}

When $m>2$, the quantities $L_{i,j}$ for more general domains can be
represented geometrically as the moduli of images of maps of $U$ onto
various vertical slit rectangles, see \fref{slit-rect}.
(In the remainder of this subsection $i=\sqrt{-1}$.)  We
explain here how to compute the $L_{j,k}$'s when the domain $U$ is the
upper half-plane.
Consider the nodes formed from the intervals $[x_1,x_2], [x_3,x_4],
\dots, [x_{2m-1},x_{2m}]$, where $x_1<x_2<\dots<x_{2m}$.
Consider the analytic function $f$ defined in the upper half-plane by
$$f(z)=\int_0^z\frac{\prod_{\ell=1}^{m-2}(w-\alpha_\ell)\,dw}{i
  \prod_{\ell=1}^{2m} (w-x_\ell)^{1/2}},$$
where
$\alpha_1,\dots,\alpha_{m-2}$ are parameters.  By the
Schwarz-Christoffel formula, if the $\alpha$'s are chosen judiciously,
then $f$ will map the upper half-plane to a rectangle with vertical
slits, such that one of the intervals $[x_{2j-1},x_{2j}]$ ($j\neq 1$)
is mapped to the top of the rectangle, the remaining intervals are
mapped to the bottom of the rectangle, and the spaces between the
intervals are mapped to the two sides of the rectangle and the
vertical slits on the bottom of the rectangle.  (The example in
\fref{slit-rect} is for $j=2$.  The case $j=1$ is similar,
except that ``top'' and ``bottom'' are reversed.)  The quantity
$L_{j,k}$ is then given by length of the image of $[x_{2k-1},x_{2k}]$
divided by the height of the rectangle:
$$L_{j,k} = i \frac{f(x_{2k})-f(x_{2k-1})}{f(x_{2j-1})-f(x_{2j-2})}.$$
(Here the choice of the $\alpha$'s, and hence $f$, depends upon $j$.)

Define the analytic function $g_p$ in the upper half-plane by
$$g_p(z) = \int_{0}^z \frac{w^p dw}{i\prod_{j=1}^{2m} (w-x_j)^{1/2}}.$$
Each of the above analytic functions $f=f_j$ is a linear combination
of $g_0,g_1,\dots,g_{m-2}$.
Let us define two $m\times (m-1)$ matrices $V'$ and $I'$ by
$$
V'_{k,p} = \Im(g_p(x_{2})-g_p(x_{2k-1}))
\ \ \ \ \text{and} \ \ \ \ 
I'_{k,p} = g_p(x_{2j})-g_p(x_{2j-1})
$$
Since for any $j\neq 1$, we can express $f_j$ as a linear
combination of the $g_p$'s, this same linear combination of the
columns of $V'$ will be nonzero (the height of the slit rectangle) at
row $j$ and zero elsewhere, and the same linear combination of the
columns of $I'$ gives the lengths of the images of the
intervals in the slit rectangle.
Thus $V'$ has rank $m-1$, and since the first row is all $0$, we may
adjoin an $m$\th column which is all $1$ to get a nonsingular matrix $V$,
and adjoin an $m$\th column to $I'$ which is all $0$ to get a matrix $I$.
Then for any vector $\vec v$, we have $L V \vec v = I \vec v$, which allows
us to compute $L$ by
$$L = I V^{-1}.$$
Thus we may compute the limiting connection probabilities completely
mechanically.
\old{
In fact there is one extra relation satisfied by $f$ and $g$,
the Riemann period relation for
$U\smallsetminus ([x_1,x_2]\cup[x_3,x_4]\cup[x_5,x_6])$.
}

\subsection{Double-dimer crossings \texorpdfstring{---}{-} multichordal \texorpdfstring{$\SLE_4$}{SLE(4)}?}

There is a lot of evidence to support the hypothesis that the scaling limit
of chains in the double-dimer model is given by $\SLE_4$, but this has not yet
been rigorously proved.  The formulas given here, which are for
double-dimer path scaling limits, agree with the formulas for $\SLE_4$
(see also \sref{GFF}), and so provide further evidence
supporting this hypothesis.

As in \sref{sle2}, we approximate the upper half-plane with
the upper half cartesian lattice $\eps\Z\times\eps\N$, and consider
a set of black nodes near the points $x_1,x_3,\dots,x_{2m-1}$ and white nodes
near the points $x_2,x_4,\dots,x_{2m}$, letting node
$i$ be the nearest lattice point $x_i^{(\eps)}\in \eps\Z\times\eps\N$
to $x_i$ which has the appropriate color.
From~\cite[\S~5.2]{MR1782431} we have
\begin{equation} \label{Xij}
\X_{i,j} = (1+o(1)) \frac{2}{\pi} \frac{\eps}{|x_i-x_j|}.
\end{equation}

For the double-dimer model we can easily recover the cross-ratio
formula that is known for bichordal $\SLE_4$
\cite[\S~8.2 and~8.3]{MR2187598} \cite{schramm-wilson} \cite[\S~4.1]{MR2253875} as follows:
\begin{align*}
 \hPr(12|34)=&\X_{1,2} \X_{3,4} \\
 \hPr(14|23)=&\X_{1,4} \X_{3,2}
\end{align*}
When $x_1<x_2<x_3<x_4$, in the scaling limit $\eps\to0$ we get
\begin{align*}
 \Pr(14|23)
 = \frac{\hPr(14|23)}{\hPr(12|34)+\hPr(14|23)}
 &\to \frac{\frac{1}{|x_1-x_4|}\frac{1}{|x_2-x_3|}}{\frac{1}{|x_1-x_2|}\frac{1}{|x_4-x_3|} + \frac{1}{|x_1-x_4|}\frac{1}{|x_2-x_3|}} \\
 &= \frac{(x_2-x_1)(x_4-x_3)}{(x_3-x_1)(x_4-x_2)}
\end{align*}
which is the formula for the cross-ratio.

For the trichordal case the relevant polynomials are
\begin{align*}
 \hPr(12|34|56)=&\X_{1,2} \X_{3,4} \X_{5,6}+\X_{1,4} \X_{3,6} \X_{5,2} \\
 \hPr(16|23|45)=&\X_{1,6} \X_{3,2} \X_{5,4}+\X_{1,4} \X_{3,6} \X_{5,2} \\
 \hPr(14|23|56)=&\X_{1,4} \X_{3,2} \X_{5,6}-\X_{1,4} \X_{3,6} \X_{5,2} \\
 \hPr(25|16|34)=&\X_{1,6} \X_{3,4} \X_{5,2}-\X_{1,4} \X_{3,6} \X_{5,2} \\
 \hPr(36|12|45)=&\X_{1,2} \X_{3,6} \X_{5,4}-\X_{1,4} \X_{3,6} \X_{5,2}.
\end{align*}

\subsection{Gaussian free-field contour lines and multichordal \texorpdfstring{$\SLE_4$}{SLE(4)}}
\label{GFF}

Schramm and Sheffield studied the contour lines of the discrete
Gaussian free field on a planar domain with boundary conditions
$+\lambda$ on one segment of the boundary and $-\lambda$ on the rest
of the boundary \cite{math.PR/0605337}, and they showed that for a
certain value of $\lambda$ the contour line of $0$ height connecting
the two special boundary points is in the scaling limit given by
chordal $\SLE_4$.  In general there can be multiple boundary points
between which the boundary heights alternate between $+\lambda$ and
$-\lambda$, and the contour lines will form a random pairing of these
boundary points.  When there are multiple boundary points, at which
the height alternates between $+\lambda$ and $-\lambda$, each of these
boundary points induces a drift term on the Brownian motion driving
the $\SLE_4$ process \cite{math.PR/0605337} \cite{schramm-sheffield:personal}.
Given the conjecture
that the scaling limit of the double-dimer model yields $\SLE_4$, and
the fact that we know the scaling limit of the connection
probabilities in the double-dimer model with multiple strands, it is
natural to compare these probabilities to the connection probabilities
of the contour lines in the discrete Gaussian free field.
We already have a natural guess for
these contour line connection probabilities (namely that they equal
the connection probabilities in the double-dimer model), and it is not
difficult to verify that this guess is in fact correct using the SLE
formulation.  In this subsection we prove
\begin{theorem} \label{GFF-pairing}
  Consider the scaling limit of the discrete Gaussian free field on a
  planar lattice in the upper half-plane, with boundary heights that
  alternate between $+\lambda$ and $-\lambda$ at the points
  $x_1<\dots<x_{2n}$, for the value of $\lambda$ determined in
  \cite{math.PR/0605337}.  Then the probability distribution of the
  manner in which the height-$0$ contour lines connect up with one
  another coincides with the connection probability distribution for
  the scaling limit of the double-dimer model in the upper half-plane
  with boundary nodes that alternate between black and white at
  locations $x_1,\dots,x_{2n}$.
\end{theorem}

Recall from \sref{ddimer-poly} and \sref{ddimer-integral} that the
probability of any double-dimer pairing is a linear combination of
$D_S/D_\varnothing$'s, where $S$ ranges over balanced subsets of nodes.
(An expression for $D_S$ is recalled below.)  We therefore start with
a lemma giving the scaling limit of $D_S/D_\varnothing$.
\begin{lemma}
  In the scaling limit of the double-dimer model in the upper
  half-plane, if the nodes are located at $x_1<\dots<x_{2n}$ and the
  node colors alternate between black and white, and if $S$ is a
  balanced subset of the nodes, then
$$\frac{D_S}{D_\varnothing} = \prod_{i\in S, j\notin S}
|x_j-x_i|^{(-1)^{1+i+j}}.$$
\end{lemma}

\begin{proof}
  Recall that for a set $S$ of nodes with equal numbers of odd and
  even nodes, we defined
\begin{align*}
 D_S &=
\det [(1_{i,j\in S} + 1_{i,j\notin S})\times (-1)^{(|i-j|-1)/2}X_{i,j}]^{i=1,3,\dots 2n-1}_{j=2,4,\dots 2n} \\
\intertext{which, in the scaling limit, using \eref{Xij}, when $x_1<x_2<\dots<x_{2n}$, is given by}
&=
\det \left[(1_{i,j\in S} + 1_{i,j\notin S})\times \frac{(-1)^{(|i-j|-1)/2}}{x_j-x_i}\right]^{i=1,3,\dots 2n-1}_{j=2,4,\dots 2n} \\
&=
\pm
\det \left[\frac{1}{x_j-x_i}\right]^{\text{$i$ odd, $i\in S$}}_{\text{$j$ even, $j\in S$}}
\det \left[\frac{1}{x_j-x_i}\right]^{\text{$i$ odd, $i\notin S$}}_{\text{$j$ even, $j\notin S$}}
\intertext{which is a product of two Cauchy determinants}
&= \pm \frac{\displaystyle\prod_{\substack{i<j\\ \text{$i+j$ even} \\ \text{$i,j\in S$ or $i,j\notin S$}}}(x_j-x_i) }{\displaystyle\prod_{\substack{i<j\\\text{$i+j$ odd}\\\text{$i,j\in S$ or $i,j\notin S$}}}(x_j-x_i)}
\end{align*}
where the $i,j\in S$ factors are from the first determinant, and the
$i,j\notin S$ factors are from the second determinant.  Since each
multiplicand is positive and by \lref{Sdet} $D_S$ is a ratio of
positive quantities, the global sign is $+$.

We may evaluate the ratio
$$
\frac{D_S}{D_{\varnothing}} = \prod_{i<j} (x_j-x_i)^{\alpha_{i,j}}$$
where
$$\alpha_{i,j} = \begin{cases} 0, & \text{if $i,j\in S$ or $i,j\notin S$,}
  \\ (-1)^{1+i+j}, & \text{if $i\in S$ and $j\notin S$ or vice versa.}
\end{cases}$$
\end{proof}
We remark that products of the form $\prod_{i<j}(x_j-x_i)^{*_{i,j}}$
also arise in \cite{MR2253875}.

\begin{proof}[Proof of \tref{GFF-pairing}]
  The strategy is to show that each $D_S/D_\varnothing$ is a martingale
  for the multichordal $\SLE_4$ diffusion associated with the contour
  lines of the discrete Gaussian free field \cite{math.PR/0605337}
  \cite{schramm-sheffield:personal}.
  Then each of the double-dimer crossing probabilities will be a
  martingale for the diffusion, and since each of them has the right
  boundary conditions, this will show that the contour line crossing
  probabilities are equal to the double-dimer crossing probabilities.
  By symmetry considerations, it suffices to consider the $\SLE_4$
  process started at $x_1$ where the other $x_i$'s are force points.
To check that $M$ is a local martingale we need to verify that
\begin{align*}
0 \stackrel{?}{=} \partial_t M &=
  \sum_{i>1} \frac{2}{x_i-x_1} \frac{\partial M}{\partial x_i}
  + \frac{\partial M}{\partial x_1} \sum_{i>1} \frac{(-1)^i 2}{x_i-x_1}
  + \frac{4}{2}\frac{\partial^2 M}{\partial x_1 ^2}. \\
  \intertext{The terms in the first sum come from the Loewner evolution
  of the points $x_i$ ($i\neq1$) when the SLE$_4$ at $x_1$ is developed,
  the terms in the second sum arise from the drifts that the points at
  $x_i$ ($i\neq1$) induce on the Brownian motion driving the SLE$_4$ at $x_1$,
  and the last term comes from the quadratic variation of the Brownian
  motion driving the SLE$_4$ process at $x_1$.
 Upon dividing by $2M$ the equation becomes}
0 \stackrel{?}{=}\frac{\partial_t M}{2 M} &=
  \sum_{i>1} \frac{1}{x_i-x_1} \frac{\partial \log M}{\partial x_i}
  + \frac{\partial \log M}{\partial x_1} \sum_{i>1} \frac{(-1)^i }{x_i-x_1}
  + \frac{\partial^2 \log M}{\partial x_1 ^2}
  + \Big(\frac{\partial \log M}{\partial x_1}\Big)^2.
\end{align*}
If $M=D_S/D_\varnothing$, then $\log M =
\sum_{i<j} \alpha_{i,j} \log(x_j-x_i)$, so the first term is
$$ \sum_{j<k} \frac{\alpha_{j,k}}{(x_k-x_j)(x_k-x_1)} +\! \sum_{1<j<k} \frac{-\alpha_{j,k}}{(x_k-x_j)(x_j-x_1)} = \sum_{1<j<k} \frac{-\alpha_{j,k}}{(x_j-x_1)(x_k-x_1)} +\! \sum_{1<k} \frac{\alpha_{1,k}}{(x_k-x_1)^2}, $$
the second term is
$$ \left(-\sum_{1<j} \frac{\alpha_{1,j}}{x_j-x_1}\right)\left(\sum_{1<k}\frac{(-1)^k}{x_k-x_1}\right), $$
the third term is
$$ \sum_{1<j} \frac{-\alpha_{1,j}}{(x_j-x_1)^2},$$
and the fourth term is
$$ \left(\sum_{1<j} \frac{-\alpha_{1,j}}{x_j-x_1}\right)^2.$$
Adding these up we get
$$
\sum_{1<j} \frac{\alpha_{1,j}-\alpha_{1,j}(-1)^j-\alpha_{1,j}+\alpha_{1,j}^2}{(x_j-x_1)^2}
+
\sum_{1<j<k} \frac{-\alpha_{j,k}-\alpha_{1,j}(-1)^k - \alpha_{1,k}(-1)^j + 2\alpha_{1,j}\alpha_{1,k}}{(x_j-x_1)(x_k-x_1)}.
$$
If $\alpha_{1,j}\neq 0$, then $\alpha_{1,j}=(-1)^j$, so the
numerators in the first sum are all $0$.  For the second sum, there
are a few cases to consider.  If $\alpha_{j,k}\neq 0$ then
$\alpha_{j,k}=(-1)^{1+j+k}$, and either $\alpha_{1,j}=0$ or else
$\alpha_{1,k}=0$ but not both, so the numerator simplifies to $0$.  If
$\alpha_{j,k}=0$ then either $\alpha_{1,j}=\alpha_{1,k}=0$ (in which
case the numerator simplifies to $0$), or else $\alpha_{1,j}=(-1)^{j}$
and $\alpha_{1,k}=(-1)^{k}$, and the numerator again simplifies to
$0$.  Thus $D_S/D_\varnothing$ is a local martingale.  Since $D_S/D_\varnothing$
is the probability that no node in $S$ is paired to a node in $S^c$
(in the double-dimer model),
it is bounded, and therefore it is actually a martingale.
\end{proof}

\section*{Acknowledgement}
We thank Dylan Thurston for helpful discussions.

\appendix

\section{The response matrix}\label{response}

To a pair $(\G,\No)$ consisting of a graph and a subset $\No$ of its
vertices is associated an $|\No|\times|\No|$ matrix $\Lambda$, the
\emph{response\/} matrix, or \emph{Dirichlet-to-Neumann\/} matrix.  It
is a real symmetric matrix defined from the Laplacian matrix as
follows.  The variables $L_{i,j}$ used in this article are negatives
of the entries of this matrix: $L_{i,j}=-\Lambda_{i,j}$.

\subsection{Definition}

Let $\Delta$ be the Laplacian matrix of $\G$, defined on functions on
the vertices of $\G$ by $\Delta(f)(v) = \sum_{w\sim
  v}c_{v,w}(f(v)-f(w))$, where the sum is over neighbors $w$ of $v$
and $c_{v,w}$ is the edge weight.  Let $\No$ be the set of nodes of
$\G$ and $\I=\V\smallsetminus \No$ the inner vertices.  Then $\R^\V=\R^\No\oplus\R^{\I}$.  In this
ordering, we can write
$$\Delta=\begin{pmatrix}F&G\\H&K\end{pmatrix}$$
where $K$ is an $|\I|\times|\I|$ matrix.

Let $V\in\R^{\No}$ be a choice of potentials at the nodes.
If we hold the nodes at these potentials, the resulting current flow
in the network will have divergence zero except at the nodes,
that is,
$$\begin{pmatrix}F&G\\H&K\end{pmatrix}
\begin{pmatrix}v_\No\\ v_{\I}\end{pmatrix}=
\begin{pmatrix}c_\No\\0\end{pmatrix},$$
where $c_\No\in\R^\No$
is the vector of currents exiting the network at the nodes.

In particular we must have $v_\I=-K^{-1}Hv_\No$ and
so $(F-GK^{-1}H)v_\No=c_\No.$
(The matrix~$K$ is invertible since it is the Dirichlet Laplacian on $\I$
with nonempty boundary $\No$.)
The matrix $\Lambda=F-GK^{-1}H$ satisfies $$\Lambda v_\No = c_\No$$ and is called the
\emph{response matrix,} or \emph{Dirichlet-to-Neumann matrix\/} of the pair $(\G,\No)$,
and is also known as the Schur complement $\Delta/K$.
It is a symmetric positive semidefinite matrix from
$\R^\No$ to $\R^\No$ with kernel consisting
of the ``all ones'' vector \cite{\CdV}.

The following lemma is a special case of a more general determinant formula,
but we need this special case in \sref{groves}, so we prove it here.
\begin{lemma}\label{putree}
  $\Pr(\tree)/\Pr(\unc) = \det \tilde\Lambda,$ where $\tilde\Lambda$ is obtained from
  $\Lambda$ by removing any row and column.
\end{lemma}
\begin{proof}
  Since $K$ is the Dirichlet Laplacian for the graph $\G/\No$ in which
  all vertices in $\No$ have been identified, by the matrix-tree
  theorem (often attributed to Kirchhoff, see
  \cite[Chapter~34]{MR1207813}), $\det K$ is
  the weighted sum of spanning trees of $\G/\No$.  These lift to
  completely disconnected groves of $\G$, so $\det K = Z(\unc)$, where
  $Z(\tau)$ denotes the weighted sum of groves of partition type~$\tau$,
  in this case $1|2|\cdots|n$.

  Let $\Delta_\eps=\Delta+\eps I_\No$; recall that this is the
  Dirichlet Laplacian for the graph $\G'$ obtained from $\G$ by
  adjoining an external vertex connected to every vertex in $\No$ by
  an edge of weight~$\eps$, i.e., $\Delta_\eps$ is obtained from
  $\Delta(\G')$ by striking out the row and column corresponding to
  the adjoined external vertex, so that $\lim_{\eps\to0}
  \det(\Delta+\eps I_\No)/\eps = n Z(\tree)$.

  Let $\Lambda_\eps=\Delta_\eps/K = F+\eps I_\No-G K^{-1}H$.  Then
  $$\Delta_\eps = \begin{pmatrix}\Lambda_\eps&G K^{-1}\\0&I_\I\end{pmatrix}
  \begin{pmatrix}I_\No&0\\H&K\end{pmatrix},$$
  so $\det\Delta_\eps=\det \Lambda_\eps \det K$ (a standard fact about Schur complements).

  But if $\tilde \Lambda$ denotes $\Lambda$ with a row and column removed, then
  $\det\tilde\Lambda = \frac1n \lim_{\eps\to0} (\det\Lambda_\eps)/\eps
  = \frac1n \lim_{\eps\to0} (\det\Delta_\eps)/\eps / \det K= Z(\tree)/Z(\unc)$.
\end{proof}

\subsection{Relation between the \texorpdfstring{$L_{i,j}$}{L}'s and \texorpdfstring{$R_{i,j}$}{R}'s}

The $\Lambda$ matrix is related to the matrix of pairwise resistances
between the nodes as follows.  We assume here for convenience that the
electrical network is connected.
\begin{proposition} \label{LtoR}
We have $$R_{i,j}= (\delta_i-\delta_j)^T \Lambda^{-1}(\delta_i-\delta_j).$$
\end{proposition}

Note that the right hand side is well defined since
$\delta_i-\delta_j$ is in the image of $\Lambda$, and the inverse image of
this vector is well-defined up to addition of an element of the kernel
of $\Lambda$, which is perpendicular to $\delta_i-\delta_j$.

\begin{proof}
  Find $v\in\R^V$ such that $\Lambda v=-\delta_i+\delta_j$.  This represents
  the (unique up to an additive constant) choice of potentials for
  which there is one unit of current flowing into the circuit at $i$
  and out at $j$.  By definition, the resistance between $i$ and $j$
  is the difference in potentials at $i$ and $j$ for this current flow.
\end{proof}

The equation of \pref{LtoR} can be given in a
computationally simpler form as follows.  When finding $v\in\R^V$ such
that $\Lambda v=-\delta_i+\delta_j$, we may choose $v$ so that
$v_n=0$, and let $\tilde v$ denote $v$ with $v_n$ dropped.  Let
$\tilde \Lambda$ be the matrix obtained from $\Lambda$ by deleting the
last row and column.  Then $\tilde\Lambda\tilde v =
-\delta_i+\delta_j$, where $\delta_n$ is interpreted as $0$.  By
the matrix-tree theorem, $\tilde\Lambda$ is invertible since
$\det\tilde\Lambda$ is the weighted sum of spanning trees with edge
weights given by $-\Lambda_{i,j}$.  Thus
$R_{i,n}=\tilde \Lambda^{-1}_{i,i},$ and if $i,j\neq n$ then
\begin{equation}\label{ltor}
  R_{i,j}= (\delta_i-\delta_j)^T \tilde \Lambda^{-1}(\delta_i-\delta_j).
\end{equation}

Reciprocally to \pref{LtoR}
we can write $\Lambda$ in terms of $R$ as follows:

\begin{proposition}
We have $\tilde \Lambda^{-1}_{i,i}=R_{i,n},$ and
$$\tilde \Lambda^{-1}_{i,j}=\frac12(R_{i,n}+R_{j,n}-R_{i,j}).$$
\end{proposition}

\begin{proof}
This is just the inverse of the transformation in
\eqref{ltor}.
\end{proof}

\section{Complete list of grove probabilities on \texorpdfstring{$4$}{4} nodes}\label{4nodes}

Here we give the formulas for the planar partition probabilities when
there are $4$ nodes.  There are $14$ planar partitions on $4$ nodes, but for
each class of planar partitions that may be obtained from one another
by cyclically rotating and/or reversing the indices, we need only give
the formula for one representative planar partition in the class.

\begin{align*}
\pt{1234} =& 1,\\
\pt{2|143} =& \textstyle\frac12 R_{2,3} + \frac12 R_{1,2} - \frac12 R_{1,3},\\
\pt{14|23}=& \textstyle\frac12 R_{1,3} + \frac12 R_{2,4} - \frac12 R_{1,4} - \frac12 R_{2,3},\\
\pt{1|2|34}=&(
+R_{1,2}R_{1,3}
+R_{1,2}R_{2,3}
+R_{1,2}R_{1,4}
+R_{1,2}R_{2,4}
-2 R_{1,2}R_{3,4}
-R_{1,2}^2\\&
+R_{1,3}R_{2,4}
+R_{1,4}R_{2,3}
-R_{1,3}R_{1,4}
-R_{2,3}R_{2,4})/4,
\\
\pt{1|3|24}=&(
+ R_{1,3} R_{1,2}
+ R_{1,3} R_{2,3}
+ R_{1,3} R_{1,4}
+ R_{1,3} R_{3,4}
- 2 R_{1,3} R_{2,4}
- R_{1,3}^2
\\&
+ R_{1,2} R_{3,4}
+ R_{1,4} R_{2,3}
- R_{1,2} R_{1,4}
- R_{2,3} R_{3,4})/4,
\\
\pt{1|2|3|4}=&(
+ R_{1,2} R_{2,3} R_{3,4}
+ R_{1,2} R_{2,4} R_{4,3}
+ R_{1,3} R_{3,2} R_{2,4}
+ R_{1,3} R_{3,4} R_{4,2}\\&
+ R_{1,4} R_{4,2} R_{2,3}
+ R_{1,4} R_{4,3} R_{3,2}
+ R_{2,1} R_{1,3} R_{3,4}
+ R_{2,1} R_{1,4} R_{4,3}\\&
+ R_{2,3} R_{3,1} R_{1,4}
+ R_{2,4} R_{4,1} R_{1,3}
+ R_{3,1} R_{1,2} R_{2,4}
+ R_{3,2} R_{2,1} R_{1,4}
\\&
- R_{1,2} R_{2,3} R_{3,1}
- R_{1,2} R_{2,4} R_{4,1}
- R_{1,3} R_{3,4} R_{4,1}
- R_{2,3} R_{3,4} R_{4,2}
\\&
- R_{1,2} R_{3,4}^2
- R_{1,2}^2 R_{3,4}
- R_{1,3} R_{2,4}^2
- R_{1,3}^2 R_{2,4}
- R_{1,4} R_{2,3}^2
- R_{1,4}^2 R_{2,3})/4.
\end{align*}
Equivalently,
\begin{alignat*}{2}
\pu{1234} =& \KL_{1234} &=&
 L_{1,2} L_{1,3} L_{1,4} +
 L_{1,2} L_{2,3} L_{2,4} +
 L_{1,3} L_{2,3} L_{3,4} +
 L_{1,4} L_{2,4} L_{3,4} +
\\&&&
 L_{1,2} L_{2,3} L_{3,4} +
 L_{1,3} L_{2,3} L_{2,4} +
 L_{1,2} L_{2,4} L_{3,4} +
 L_{1,4} L_{2,3} L_{2,4} +
\\&&&
 L_{1,2} L_{1,3} L_{2,4} +
 L_{1,2} L_{1,4} L_{2,3} +
 L_{1,3} L_{2,4} L_{3,4} +
 L_{1,4} L_{2,3} L_{3,4} +
\\&&&
 L_{1,3} L_{1,4} L_{2,3} +
 L_{1,2} L_{1,3} L_{3,4} +
 L_{1,3} L_{1,4} L_{2,4} +
 L_{1,2} L_{1,4} L_{3,4},\\
\pu{2|143} =&\KL_{2|143} + \KL_{13|24}&\ =&
 L_{1,3} L_{1,4} + L_{1,4} L_{3,4} + L_{1,3} L_{3,4} + L_{1,3} L_{2,4},\\
\pu{14|23}=&\KL_{14|23} - \KL_{13|24}&=&
 L_{1,4} L_{2,3} - L_{1,3} L_{2,4},\\
\pu{1|2|34}=&\KL_{1|2|34}&=&L_{3,4},\\
\pu{1|3|24}=&\KL_{1|3|24}&=&L_{2,4},\\
\pu{1|2|3|4}=&\KL_{1|2|3|4}&=&1.
\end{alignat*}

\pdfbookmark[1]{References}{bib}
\bibliographystyle{hmralpha}
\bibliography{bc}

\begin{thebibliography}{DFGG97}

\bibitem[ASA02]{MR1928682}
Louis-Pierre Arguin and Yvan Saint-Aubin.
\newblock Non-unitary observables in the 2d critical {I}sing model.
\newblock {\em Physics Letters B}, 541(3-4):384--389, 2002,
  \arXiv{hep-th/0109138}. \MR{MR1928682 (2003h:82016)}

\bibitem[BBK05]{MR2187598}
Michel Bauer, Denis Bernard, and Kalle Kyt{\"o}l{\"a}.
\newblock Multiple {S}chramm-{L}oewner evolutions and statistical mechanics
  martingales.
\newblock {\em J. Stat.\ Phys.}, 120(5-6):1125--1163, 2005,
  \arXiv{math-ph/0503024}. \MR{MR2187598 (2007d:82032)}

\bibitem[BSW07]{math.CA/0703826}
Nathaniel~D. Blair-Stahn and David~B. Wilson.
\newblock {The electrical response matrix of a regular {$2n$}-gon}, 2007,
  \arXiv{math.CA/0703826}.
\newblock \textit{Proc.\ Amer.\ Math.\ Soc.}, to appear.

\bibitem[Car92]{MR92m:82048}
John~L. Cardy.
\newblock Critical percolation in finite geometries.
\newblock {\em J. Phys.\ A}, 25(4):L201--L206, 1992, \arXiv{math-ph/9910002}.
  \MR{MR1151081 (92m:82048)}

\bibitem[Car07]{cardy:ade}
John Cardy.
\newblock {ADE} and {SLE}.
\newblock {\em J. Phys.\ A}, 40(7):1427--1438, 2007, \arXiv{math-ph/0610030}.
  \MR{MR2303259 (2008k:82023)}

\bibitem[CdV98]{MR1652692}
Yves Colin~de Verdi{\`e}re.
\newblock {\em Spectres de Graphes}, volume~4 of {\em Cours Sp\'ecialis\'es
  [Specialized Courses]}.
\newblock Soci\'et\'e Math\'ematique de France, Paris, 1998. \MR{MR1652692
  (99k:05108)}

\bibitem[CIM98]{MR1657214}
E.~B. Curtis, D.~Ingerman, and J.~A. Morrow.
\newblock Circular planar graphs and resistor networks.
\newblock {\em Linear Algebra Appl.}, 283(1-3):115--150, 1998. \MR{MR1657214
  (99k:05096)}

\bibitem[Ciu97]{ciucu}
Mihai Ciucu.
\newblock Enumeration of perfect matchings in graphs with reflective symmetry.
\newblock {\em J. Combin.\ Theory Ser.\ A}, 77(1):67--97, 1997. \MR{MR1426739
  (98a:05112)}

\bibitem[CN07]{camia-newman:proof-converge}
Federico Camia and Charles~M. Newman.
\newblock Critical percolation exploration path and {${\rm SLE}\sb 6$}: a proof
  of convergence.
\newblock {\em Probab.\ Theory Related Fields}, 139(3-4):473--519, 2007,
  \arXiv{math/0604487}. \MR{MR2322705 (2008k:82040)}

\bibitem[CS04]{MR2097339}
Gabriel~D. Carroll and David Speyer.
\newblock The cube recurrence.
\newblock {\em Electron.\ J. Combin.}, 11(1):Research Paper 73, 31 pp., 2004,
  \arXiv{math.CO/0403417}. \MR{MR2097339 (2005f:05007)}

\bibitem[DFGG97]{MR1462755}
P.~Di~Francesco, O.~Golinelli, and E.~Guitter.
\newblock Meanders and the {T}emperley-{L}ieb algebra.
\newblock {\em Comm.\ Math.\ Phys.}, 186(1):1--59, 1997,
  \arXiv{hep-th/9602025}. \MR{MR1462755 (99f:82028)}

\bibitem[DS84]{MR920811}
Peter~G. Doyle and J.~Laurie Snell.
\newblock {\em Random Walks and Electric Networks}, volume~22 of {\em Carus
  Mathematical Monographs}.
\newblock Mathematical Association of America, Washington, DC, 1984,
  \arXiv{math.PR/0001057}. \MR{MR920811 (89a:94023)}

\bibitem[DS87]{duplantier-saleur:dense-saw}
B.~Duplantier and H.~Saleur.
\newblock Exact critical properties of two-dimensional dense self-avoiding
  walks.
\newblock {\em Nuclear Phys.\ B}, 290(3):291--326, 1987. \MR{MR919521
  (89g:82049)}

\bibitem[Dub06]{MR2253875}
Julien Dub{\'e}dat.
\newblock Euler integrals for commuting {SLE}s.
\newblock {\em J. Stat.\ Phys.}, 123(6):1183--1218, 2006,
  \arXiv{math.PR/0507276}. \MR{MR2253875 (2007g:82027)}

\bibitem[Dup87]{duplantier:On}
Bertrand Duplantier.
\newblock Critical exponents of {M}anhattan {H}amiltonian walks in two
  dimensions, from {P}otts and {${\rm O}(n)$} models.
\newblock {\em J. Statist.\ Phys.}, 49(3-4):411--431, 1987. \MR{MR926190
  (89b:82025)}

\bibitem[Dup89]{duplantier:networks}
Bertrand Duplantier.
\newblock Statistical mechanics of polymer networks of any topology.
\newblock {\em J. Statist.\ Phys.}, 54(3-4):581--680, 1989. \MR{MR988554
  (90e:82061)}

\bibitem[Dup92]{duplantier:lerw}
Bertrand Duplantier.
\newblock Loop-erased self-avoiding walks in two dimensions: exact critical
  exponents and winding numbers.
\newblock {\em Physica A}, 191:516--522, 1992.

\bibitem[Dup06]{duplantier:conformal-geometry}
Bertrand Duplantier.
\newblock Conformal random geometry.
\newblock In A.~Bovier, F.~Dunlop, A.~van Enter, F.~den Hollander, and
  J.~Dalibard, editors, {\em Mathematical Statistical Physics}, Lecture Notes
  of the Les Houches Summer School Session LXXXIII, 2005, pages 101--217.
  Elsevier, 2006, \arXiv{math-ph/0608053}.

\bibitem[Fom01]{MR1837248}
Sergey Fomin.
\newblock Loop-erased walks and total positivity.
\newblock {\em Trans.\ Amer.\ Math.\ Soc.}, 353(9):3563--3583 (electronic),
  2001, \arXiv{math.CO/0004083}. \MR{MR1837248 (2002f:15030)}

\bibitem[Jon83]{MR696688}
V.~F.~R. Jones.
\newblock Index for subfactors.
\newblock {\em Invent.\ Math.}, 72(1):1--25, 1983. \MR{MR696688 (84d:46097)}

\bibitem[Kas67]{MR0253689}
P.~W. Kasteleyn.
\newblock Graph theory and crystal physics.
\newblock In {\em Graph Theory and Theoretical Physics}, pages 43--110.
  Academic Press, London, 1967. \MR{MR0253689 (40 \#6903)}

\bibitem[Ken97]{MR1473567}
Richard Kenyon.
\newblock Local statistics of lattice dimers.
\newblock {\em Ann.\ Inst.\ H. Poincar\'e Probab.\ Statist.}, 33(5):591--618,
  1997, \arXiv{math.CO/0105054}. \MR{MR1473567 (99b:82039)}

\bibitem[Ken00]{MR1782431}
Richard Kenyon.
\newblock Conformal invariance of domino tiling.
\newblock {\em Ann.\ Probab.}, 28(2):759--795, 2000, \arXiv{math-ph/9910002}.
  \MR{MR1782431 (2002e:52022)}

\bibitem[Kir90]{Kirchhoff}
Gustav Kirchhoff.
\newblock Beweis der {E}xistenz des {P}otentials das an der {G}renze des
  betrachteten {R}aumes gegebene {W}erthe hat f\"ur den {F}all dass diese
  {G}renze eine \"uberall convexe {F}l\"ache ist.
\newblock {\em Acta Math.}, 14(1):179--183, 1890. \MR{MR1554794 {}}

\bibitem[KL07]{math.PR/0605159}
Michael~J. Kozdron and Gregory~F. Lawler.
\newblock The configurational measure on mutually avoiding {SLE} paths.
\newblock In {\em Universality and Renormalization}, volume~50 of {\em Fields
  Inst.\ Commun.}, pages 199--224. Amer.\ Math.\ Soc., Providence, RI, 2007,
  \arXiv{math.PR/0605159}. \MR{MR2310306 (2008m:60079)}

\bibitem[KS91]{MR1105701}
Ki~Hyoung Ko and Lawrence Smolinsky.
\newblock A combinatorial matrix in {$3$}-manifold theory.
\newblock {\em Pacific J. Math.}, 149(2):319--336, 1991. \MR{MR1105701
  (92d:57008)}

\bibitem[Kuo04]{kuo-condense}
Eric~H. Kuo.
\newblock Applications of graphical condensation for enumerating matchings and
  tilings.
\newblock {\em Theoretical Computer Science}, 319(1-3):29--57, 2004,
  \arXiv{math.CO/0304090}.

\bibitem[KW08]{KW-Pfaffian}
Richard~W. Kenyon and David~B. Wilson.
\newblock Combinatorics of tripartite boundary connections for trees and
  dimers, 2008, \arXiv{0811.1766}.

\bibitem[LSW04]{MR2044671}
Gregory~F. Lawler, Oded Schramm, and Wendelin Werner.
\newblock Conformal invariance of planar loop-erased random walks and uniform
  spanning trees.
\newblock {\em Ann.\ Probab.}, 32(1B):939--995, 2004, \arXiv{math.PR/0112234}.
  \MR{MR2044671 (2005f:82043)}

\bibitem[OB88]{ohno-binder:networks}
Kaoru Ohno and Kurt Binder.
\newblock Scaling theory of star polymers and general polymer networks in bulk
  and semi-infinite good solvents.
\newblock {\em J. Physique}, 49(8):1329--1351, 1988. \MR{MR966320 (89j:82053)}

\bibitem[Pem95]{MR1410532}
Robin Pemantle.
\newblock Uniform random spanning trees.
\newblock In J.~Laurie Snell, editor, {\em Topics in Contemporary Probability
  and its Applications}, Probab. Stochastics Ser., pages 1--54. CRC, Boca
  Raton, FL, 1995. \MR{MR1410532 (97i:60013)}

\bibitem[PS05]{MR2144860}
T.~Kyle Petersen and David Speyer.
\newblock An arctic circle theorem for {G}roves.
\newblock {\em J. Combin.\ Theory Ser.\ A}, 111(1):137--164, 2005,
  \arXiv{math.CO/0407171}. \MR{MR2144860 (2006m:05055)}

\bibitem[Sal86]{saleur:saw}
H.~Saleur.
\newblock New exact exponents for two-dimensional self-avoiding walks.
\newblock {\em J. Phys.\ A}, 19:L807--L810, 1986.

\bibitem[Sch00]{MR1776084}
Oded Schramm.
\newblock Scaling limits of loop-erased random walks and uniform spanning
  trees.
\newblock {\em Israel J. Math.}, 118:221--288, 2000, \arXiv{math/9904022}.
  \MR{MR1776084 (2001m:60227)}

\bibitem[Sch07]{MR2334202}
Oded Schramm.
\newblock Conformally invariant scaling limits (an overview and a collection of
  problems).
\newblock In {\em International {C}ongress of {M}athematicians. {V}ol.~{I}},
  pages 513--543. Eur.\ Math.\ Soc., Z\"urich, 2007, \arXiv{math/0602151}.
  \MR{MR2334202 (2008j:60237)}

\bibitem[Smi01]{MR1851632}
Stanislav Smirnov.
\newblock Critical percolation in the plane: conformal invariance, {C}ardy's
  formula, scaling limits.
\newblock {\em C. R. Acad.\ Sci.\ Paris S\'er.\ I Math.}, 333(3):239--244,
  2001. \MR{MR1851632 (2002f:60193)}

\bibitem[Smi07]{smirnov:ising}
Stanislav Smirnov.
\newblock Conformal invariance in random cluster models. {I}. {Holomorphic}
  fermions in the {Ising} model, 2007, \arXiv{0708.0039}.

\bibitem[Spi76]{MR0388547}
Frank Spitzer.
\newblock {\em Principles of Random Walk}.
\newblock Springer-Verlag, New York, second edition, 1976.
\newblock Graduate Texts in Mathematics, Vol.\ 34. \MR{MR0388547 (52 \#9383)}

\bibitem[SS06]{math.PR/0605337}
Oded Schramm and Scott Sheffield.
\newblock {Contour lines of the two-dimensional discrete Gaussian free field},
  2006, \arXiv{math.PR/0605337}.

\bibitem[SS07]{schramm-sheffield:personal}
Oded Schramm and Scott Sheffield, 2007.
\newblock In preparation.

\bibitem[Sta86]{MR847717}
Richard~P. Stanley.
\newblock {\em Enumerative Combinatorics, Volume 1}.
\newblock Wadsworth \& Brooks/Cole Advanced Books \& Software, 1986.
  \MR{MR847717 (87j:05003)}

\bibitem[Sta99]{MR1676282}
Richard~P. Stanley.
\newblock {\em Enumerative Combinatorics, Volume 2}.
\newblock Cambridge Studies in Advanced Mathematics \#62. Cambridge University
  Press, Cambridge, 1999.
\newblock With an appendix by Sergey Fomin. \MR{MR1676282 (2000k:05026)}

\bibitem[SW05]{schramm-wilson}
Oded Schramm and David~B. Wilson.
\newblock {SLE}, quadrangles, and curvilinear triangles, 2005.

\bibitem[vLW92]{MR1207813}
J.~H. van Lint and R.~M. Wilson.
\newblock {\em A Course in Combinatorics}.
\newblock Cambridge University Press, Cambridge, 1992. \MR{MR1207813
  (94g:05003)}

\end{thebibliography}

\end{document}